\documentclass[11pt]{article}
\usepackage{amsmath,amsthm,amsfonts,amscd,bezier}
\usepackage{amssymb}
\usepackage{enumerate}
\usepackage[T1]{fontenc}
\usepackage[usenames]{color}

\newtheorem{lemma}{Lemma}[section]

\newtheorem{theorem}[lemma]{Theorem}
\newtheorem{remark}[lemma]{Remark}
\newtheorem{proposition}[lemma]{Proposition}
\newtheorem{corollary}[lemma]{Corollary}
\newtheorem{definition}[lemma]{Definition}
\newtheorem{example}[lemma]{Example}
\newcommand{\Dem}{\noindent{\sc Proof:\ \ }}
\newcommand{\cqd}{{\hfill $\rule{2mm}{2mm}$}\vspace{1cm}}

\title{On the Analytic Invariants and Semiroots of Plane Branches}

\author{Marcelo Osnar Rodrigues de Abreu  \\ and \\  Marcelo Escudeiro Hernandes\thanks{The first author was partially supported by CAPES and the second one by CNPq-Brazil.} \thanks{Corresponding author: Hernandes, M. E.; email: mehernandes@uem.br }}

\begin{document}
\maketitle
\markboth{M. O. R. Abreu and M. E.
Hernandes}{On the Analytic Invariants and Semiroots of Plane Branches}

\begin{center} 2010 Mathematics Subject Classification: 14H20 (primary),
 32S10 (secondary)\end{center}

\begin{center} keywords: Plane curves, Semiroots, Analytic Invariants. \end{center}

\begin{abstract}
The value semigroup of a $k$-semiroot $C_k$ of a plane branch $C$
allows us to recover part of the value semigroup $\Gamma =\langle
v_0,\ldots ,v_g\rangle$ of $C$, that is, it is related to
topological invariants of $C$. In this paper we consider the set of
values of differentials $\Lambda_k$ of $C_k$, that is an analytical
invariant, and we show how it determine part of the set of values of
differentials $\Lambda$ of $C$. As a consequence, in a fixed
topological class, we relate the Tjurina number $\tau$ of $C$ with
the Tjurina number of $C_k$. In particular, we show that $\tau\leq
\mu-\frac{3n_g-2}{4}\mu_{g-1}$ where $n_g=gcd(v_0,\ldots ,v_{g-1})$,
$\mu$ and $\mu_{g-1}$ denote the Milnor number of $C$ and $C_{g-1}$
respectively. If $n_g=2$, we have that $\tau=\mu-\mu_{g-1}$ for any
curve in the topological class determined by $\Gamma$ that is a
generalization of a result obtained by Luengo and
Pfister.\end{abstract}

\section{Introduction}

Given $f\in\mathbb{C}\{x,y\}$ irreducible the equation $f=0$ define
a germ of irreducible plane curve, or a {\it branch}, $C_f$.

The topological class of $C_f$ can be totaly characterized by
numerical data: by the characteristic exponents of a Puiseux
parametrization, the value semigroup $\Gamma_f$, the multiplicity
sequence in the canonical resolution, for instance. All of these
numerical invariants determine each other and can be easily
computed.

The value semigroup $\Gamma_f\subseteq \mathbb{N}$ admits a conductor $\mu_f$ that is, $\mu_f-1\not\in\Gamma_f$ and $\mu_f+\mathbb{N}\subseteq\Gamma_f$. So, $\Gamma_f$ is determined by the finite set $\mathbb{N}\setminus\Gamma_f$ and $\mu_f$ is a topological invariant that coincide with the Milnor number of $f$, i.e., $\mu_f=\dim_{\mathbb{C}}\frac{\mathbb{C}\{x,y\}}{\langle f_x,f_y\rangle}$.

On the other hand, it is not so easy to compute analytical
invariants and to describe relation among them. In \cite{HefezHern},
the set $\Lambda_f$ of values of differentials was the principal
analytical invariant and ingredient considered to obtain normal
forms with respect to the analytical equivalence of plane branches.
The set $\Lambda_f$ is related with the Tjurina number
$\tau_f=\dim_{\mathbb{C}}\frac{\mathbb{C}\{x,y\}}{\langle
f,f_x,f_y\rangle}$, another analytical invariant of $C_f$, by the
formula $\tau_f=\mu_f-\sharp (\Lambda_f\setminus\Gamma_f)$.

Up to analytical equivalence any plane curve can be defined by a
Weierstrass polynomial $f\in\mathbb{C}\{x\}[y]$ and many topological
aspects of $C_f$ can be studied by a set of monic irreducible
polynomials $f_k\in \mathbb{C}\{x\}[y],\ 0\leq k< g$ with
$deg_y(f_k)=\frac{v_0}{e_k}$, called semiroots of $f$, such that
$\nu_f(f_k):=\dim_{\mathbb{C}}\frac{\mathbb{C}\{x,y\}}{\langle
f,f_k\rangle}=v_{k+1}$ where $\Gamma_f=\langle v_0,v_1,\ldots
,v_g\rangle$,
$deg_y(f)=v_0=\dim_{\mathbb{C}}\frac{\mathbb{C}\{x,y\}}{\langle
f,x\rangle}$ and $e_i=gcd(v_0,\ldots ,v_i)$ (see Section 2).

Several authors have considered semiroots in different contexts. For
instance, Zariski in \cite{zariskibook} considered semiroots given
by minimal polynomial of truncation of Puiseux parametrization of
$C_f$ to compute the minimal set of generators of the semigroup.
Abhyankar in \cite{abhyankar} presented an irreducibility criterion
for elements in $\mathbb{C}\{x\}[y]$ using approximate roots that
are particular cases of semiroots. In \cite{javorski}, Yavorski gave
a method to construct a miniversal family of polynomials in
$\mathbb{C}[x,y]$ that define branches with a fixed semigroup
considering semiroots. Popescu in \cite{popescu} presented several
local properties of semiroots including a proof of the Embedding
Line Theorem as an application.

It is possible to relate the topological data of $C_f$ with the
corresponding data of the branch $C_k$ associated to $f_k$. For
example, we have $$\Gamma_f=e_k\Gamma_k+\langle v_{k+1},\ldots
,v_g\rangle\ \ \ \mbox{and}\ \ \
\mu_f-1=e_k(\mu_k-1)+\sum_{j=k+1}^{g}\left (n_j -1\right )v_j,$$
where $n_j=\frac{e_{j-1}}{e_j}$, $\Gamma_k$ is the value semigroup
and $\mu_k$ is the Milnor number of $C_k$.

The aim of this work is to study relations among the analytical
invariants $\Lambda_f$, $\tau_f$, $\Lambda_k$ and $\tau_k$
associated to $C_f$ and $C_k$ respectively. For this purpose, in
Section 3 we explore the logarithmic differential forms along a
curve and we revisit results of Brian\c con, Maisonobe and Torrelli
(see \cite{briancon}) about some magic matrix associated to a plane
curve and its Jacobian ideal.

In the Section 4 we describe a common subset of
$\Lambda\setminus\Gamma$ shared by any plane branch $C$ that admits
a fixed semiroot $f_k$ and consequently we can relate the Tjurina
number $\tau$ of $C$ with the Tjurina number $\tau_k$ of $C_k$. More
explicitly, for $0\leq k<g$ we get
$\tau\leq\mu-\mu_k-e_{k+1}(n_{k+1}-2)\tau_k$.

Dimca and Greuel in \cite{dimca} conjectured that
$\frac{3}{4}\mu\leq \tau$ for any plane curve. This question
motivated several works, see \cite{alberich}, \cite{genzmer},
\cite{japa} for a proof in the irreducible case and \cite{almiron}
for any plane curve. So, for $k=g-1$ we get $\tau\leq
\mu-\mu_{g-1}-(n_g-2)\tau_{g-1}\leq
\mu-(\frac{3n_g-2}{4})\mu_{g-1}$, that is an upper bound for $\tau$
for any plane branch in a fixed topological class.

Considering the topological class given by $\Gamma=\langle
v_0,\ldots ,v_g\rangle$ with $n_g=gcd(v_0,\ldots ,v_{g-1})=2$, in
Section 5 we show that $\tau$ is constant and equal to
$\mu-\mu_{g-1}$ for any plane branch that admits the semigroup
$\Gamma$. This result is trivial if $g=1$ and it was showed by
Luengo and Pfister for $g=2$ in \cite{luengo}.

\section{Plane branches and semiroots}

In what follows we recall some classical result about plane curves using the presentation of \cite{BriesKnorrer}, \cite{abramo} and \cite{zariskibook}.

Let $\mathbb{C}\{x,y\}$ be the ring of analytic power series ring in $x$ and $y$ over $\mathbb{C}$ and $\mathcal{M}=\langle x,y\rangle$ its maximal ideal. Let $C_f\in (\mathbb{C}^2,0)$ be a singular irreducible plane curve (a plane branch) defined by $f=0$, where $f\in\mathcal{M}^n\setminus\mathcal{M}^{n+1}$ is irreducible and $n>1$ is its multiplicity.

We say that $C_f$ and $C_h$ are analytically equivalent if there
exist neighborhoods $U,V\subseteq (\mathbb{C}^2,0)$ and an analytic
isomorphism $\Psi: (\mathbb{C}^2,0)\rightarrow (\mathbb{C}^2,0)$
such that $\Psi(U\cap C_f)=V\cap C_h$. If $\Psi$ is just a
homeomorphism we say that $C_f$ and $C_h$ are topologically
equivalent.

By the Weierstrass Preparation Theorem and a change of coordinates
if it is necessary, any plane branch is analytically equivalent to a
branch defined by
$f=y^n+\sum_{i=1}^{n}c_i(x)y^{n-i}\in\mathbb{C}\{x\}[y]$ with
$n\nmid ord_x(c_n(x))=m$ and by the Newton-Puiseux theorem $C_f$
admits primitive\footnote{A parametrization $(t^n,\sum_{i}a_it^i)$
is primitive if and only if $gcd(n, i; a_i\neq 0)=1$.} Puiseux
parametrizations $(t^n,\varphi(\zeta^jt))$ with $1\leq j\leq n$
where $\varphi(t)=\sum_{i\geq m}a_it^i\in\mathbb{C}\{t\}$, $a_m\neq
0$ and $\zeta$ is a primitive $n$-th root of unity, that is,
$f(x,y)=\prod_{j=1}^{n}(y-\varphi(\zeta^j x^{\frac{1}{n}}))$.

Given a Puiseux parametrization as above we consider the following sequences of positive integers:
$$n=\beta_0=e_0,\ \ \ n_0=1;$$
$$\beta_i=\min\{k;\ a_k\neq 0\ \mbox{and}\ k\not\equiv 0\ mod\ e_{i-1}\},\ \ \ e_i=gcd(e_{i-1},\beta_i),\ \ n_i=\frac{e_{i-1}}{e_i}\ \ \ \mbox{for}\ i\geq 1.$$
As the parametrization is primitive there exists an integer $g>0$ such that $e_g=1$. The increasing finite sequence $\beta_0<\beta_1<\ldots <\beta_g$ is called the {\it characteristic sequence}.
According to Zariski (see \cite{equi}) the topological class of $C_f$ is totaly determined by its characteristic sequence.

The topological data can be done by other numerical invariants. For instance, considering the exact sequence
$$\begin{array}{ccccccccc}
\{0\} & \rightarrow & \langle f\rangle & \rightarrow & \mathbb{C}\{x,y\} & \stackrel{\varphi^*}{\rightarrow} & \mathbb{C}\{t^n,\varphi(t)\} & \rightarrow & \{0\} \\
& & & & h & \mapsto & \varphi^*(h):=h(t^n,\varphi(t)) & & \end{array}$$
we have $\mathcal{O}_f:=\frac{\mathbb{C}\{x,y\} }{\langle f\rangle}\approx\mathbb{C}\{t^n,\varphi(t)\}\subseteq\mathbb{C}\{t\}\approx \overline{\mathcal{O}}_f$. The isomorphism class of the local ring $\mathcal{O}_f$ determine the analytic class of $C_f$ and we can recover the topological data by means $\mathcal{O}_f$ considering
$$\Gamma_f=\left \{\nu_f(h):=ord_t(\varphi^*(h))=\dim_{\mathbb{C}}\frac{\mathbb{C}\{x,y\} }{\langle f,h\rangle};\ h\in \mathbb{C}\{x,y\}\setminus\langle f\rangle\right \}\subseteq\mathbb{N}.$$
In fact, $\Gamma_f$ is a numerical semigroup that admits the minimal set of generators $\{v_0,v_1,\ldots ,v_g\}$ satisfying
\begin{equation}\label{vi}v_0=\beta_0,\ \ \ v_{i+1}=n_{i}v_{i}+\beta_{i+1}-\beta_{i}=\sum_{j=1}^{i}\frac{e_{j-1}-e_j}{e_i}\beta_j+\beta_{i+1},\ \ 0\leq i< g.\end{equation}
So, $\Gamma_f$ determines and it is determined by the characteristic sequence.

The integer $n_i$ can be characterized by the property
$n_i=\min\{k>0;\ kv_i\in\langle v_0,\ldots ,v_{i-1}\rangle\}$.

\begin{remark}\label{padrao}
    Any positive integer $r$ can be uniquely expressed as $r=\sum_{i=0}^{g}s_iv_i$ with $0\leq s_i<n_i$ for $1\leq i\leq g$ and $s_0\in\mathbb{Z}$. Moreover, $r\in\Gamma_f$ if and only if, $s_0\geq 0$ (see Lemma 7.1, \cite{abramo}).
\end{remark}

As $\Gamma_f$ is a numerical semigroup it admits a conductor $\mu_f$ that is, $\mu_f-1\not\in\Gamma_f$ and $\mu_f+\mathbb{N}\subseteq\Gamma_f$. In this case, $\mu_f$ coincides with the Milnor number of $C_f$ and (see \cite{milnor} and (7.1) in \cite{abramo})
\begin{equation}\label{milnor}\mu_f=\dim_{\mathbb{C}}\frac{\mathbb{C}\{x,y\} }{J(f)}=\sum_{i=1}^{g}(n_{i}-1)v_i-v_0+1=n_gv_g-\beta_g-v_0+1=2\sharp(\mathbb{N}\setminus\Gamma),\end{equation}
where $J(f)=\langle f_x,f_y\rangle\subseteq \mathbb{C}\{x,y\}$.

Concerning the vector space $\frac{\mathbb{C}\{x,y\} }{J(f)}$, Yavorski presents the following theorem:

\begin{theorem}[Yavorski, \cite{javorski}] For any set $\mathcal{B}=\{h\in \mathbb{C}\{x\}[y];\ deg_y(h)<deg_y(f)\}$ satisfying
\begin{enumerate}
\item $\nu_f(h)\in\Gamma_f\setminus\{\mu_f-1+\gamma;\ \gamma\in\Gamma_f\}$;
\item $\nu_f(h_i)\neq\nu_f(h_j)$ for every $h_i,h_j\in\mathcal{B}$,
\end{enumerate}
we have a $\mathbb{C}$-basis $\overline{\mathcal{B}}=\left \{\overline{h}\in \frac{\mathbb{C}\{x,y\} }{J(f)};\ h\in\mathcal{B}\right \}$ for $\frac{\mathbb{C}\{x,y\} }{J(f)}$.
\end{theorem}

Notice that by the previous theorem we get \begin{equation}\label{jac-gamma}
\nu_f(J(f))-(\mu_f-1)=\Gamma_f\setminus\{0\}.\end{equation}

In \cite{zariskibook}, Zariski obtains (\ref{vi}) considering the
minimal polynomial $f_{k}(x,y)\in\mathbb{C}\{x\}[y]$ of the element
$\sum_{m\leq i<\beta_k}a_ix^{\frac{i}{n_0\cdot n_1\cdot\ldots\cdot
n_{k}}}\in\mathbb{C}(x)$ and he shows that $\nu_f(f_k)=v_{k+1}$.
This is a particular case of semiroots associated to $f$.

\begin{definition}
    Given $k\in\{0,\ldots ,g\}$ a monic polynomial $f_k\in\mathbb{C}\{x\}[y]$ of degree $n_0\cdot n_1\cdot\ldots\cdot n_{k}=\frac{n}{e_{k}}$ is called a $k$-semiroot of $f$ if $\nu_f(f_k)=v_{k+1}$ where $v_{g+1}=\infty$.
\end{definition}

Any semiroot $f_k$ is irreducible, the associated curve $C_k:=C_{f_k}$ has characteristic sequence $\frac{\beta_0}{e_{k}},\ldots ,\frac{\beta_{k}}{e_{k}}$
and it admits semigroup $\Gamma_k=\langle \frac{v_0}{e_k},\ldots ,\frac{v_k}{e_k}\rangle$ (see Theorem 5.1 in \cite{popescu}).

If $C_f$ admits a Puiseux parametrization $(t^{\beta_0},\varphi(t))$ and $f_k$ is a $k$-semiroot of $f$ then $C_k$ admits a Puiseux parametrization $\left (t^{\frac{\beta_0}{e_k}},\psi(t)\right )$ with $j^{\beta_k}(\varphi(t))=j^{\beta_k}(\psi(t^{e_k}))$ where $j^r(h)$ denotes the $r$-jet of $h$ (see \cite{popescu}, Corollary 5.3). In addition, if $f_{k}$ is a $k-$semiroot of $f$ and $f_{j}$ is a $j-$semiroot of $f_{k}$ then $f_{j}$ is a $j-$semiroot of $f$.

A set $\{f_k;\ f_k\ \mbox{is a}\ k-\mbox{semiroot of}\ f\
\mbox{with}\ 0\leq k\leq g\}$ is called a {\it complete system of
semiroots} of $f$. A complete system of semiroots can be used to
express any $h\in\mathbb{C}\{x\}[y]$ in a highlighted way as
follows.

\begin{lemma}[\cite{popescu}, Corollary 5.4]\label{decom-semirraizes} Any element
$h\in\mathbb{C}\{x\}[y]$ can be uniquely expressed as
$$
h=\sum_{finite\atop \alpha=(\alpha_{0},\ldots,\alpha_{g})\in\mathbb{N}^{g+1}}b_{\alpha}f_{0}^{\alpha_{0}}\cdots f_{g-1}^{\alpha_{g-1}}f_g^{\alpha_{g}},
$$
with $0\displaystyle \leq \alpha_{g}\leq\left[\frac{gr(h)}{gr(f)}\right], 0\leq \alpha_{k}<n_{k+1}$ for $0\leq k< g$ and $b_{\alpha}\in\mathbb{C}\{x\}$. Moreover:
\begin{enumerate}
\item $deg_y(f_{0}^{\alpha_{0}}\cdots f_g^{\alpha_{g}})\neq deg_y(f_{0}^{\gamma_{0}}\cdots f_g^{\gamma_{g}})$ if and only if $\alpha\neq \gamma$.
\item $\nu_f(b_{\alpha}f_{0}^{\alpha_{0}}\cdots f_{g-1}^{\alpha_{g-1}})\neq \nu_f(b_{\gamma}f_{0}^{\gamma_{0}}\cdots f_{g-1}^{\gamma_{g-1}})$
for $(\alpha_0,\ldots ,\alpha_{g-1})\neq(\gamma_0,\ldots
,\gamma_{g-1})$.
\end{enumerate}
\end{lemma}

\begin{remark}\label{degree} Given $h\in\mathbb{C}\{x\}[y]$ as expressed in the previous lemma with
$\alpha_g=0$ for every $\alpha\in\mathbb{N}^{g+1}$ and
$\nu_f(h)=\nu_f(b_{\alpha}f_0^{\alpha_0}\cdots
f_{g-1}^{\alpha_{g-1}} )=\sum_{i=1}^{g}\alpha_{i-1}v_i$. If
$\alpha_{i-1}=n_i-1$ for $1\leq i\leq g$ then
$$deg_y(h)= deg_y(f_0^{\alpha_0}\cdots f_{g-1}^{\alpha_{g-1}})=\sum_{i=1}^{g}(n_i-1)n_0\cdot n_1\cdot\ldots \cdot n_{i-1}=n-1.$$
\end{remark}

In what follows for $0\leq k\leq g$ we consider the multiplicative group $G_k=\{\eta\in\mathbb{C};\ \eta^{e_k}=1\}$. Remark that $\{1\}=G_g\subset\cdots\subset G_1\subset G_0$ and $n_i=\sharp\frac{G_{i-1}}{G_i}$ for all $1\leq i\leq g$.

The next result relates $\nu_f(h)$ with $\nu_k(h):=\nu_{f_k}(h)$.

\begin{proposition}\label{aval-semiraiz} Let $f_j$ be a $j-$semiroot of $f_k$ and
$f_k$ a $k-$semiroot of $f=f_g\in\mathbb{C}\{x\}[y]$ with
$\left(t^{\frac{v_0}{e_{j}}},\phi(t)\right),
\left(t^{\frac{v_0}{e_{k}}},\psi(t)\right)$ and
$(t^{v_0},\varphi(t))$ their respective parametrizations. If
$f_j\left(t^{\frac{v_0}{e_{k}}},\psi(t)\right)=at^{\frac{v_j}{e_{k}}}+\mbox{(h.o.t.)}$,
then $f_j(t^{v_0},\varphi(t))=at^{v_j}+\mbox{(h.o.t.)}$.

In particular, for any $h\in\mathbb{C}\{x\}[y]$ with $deg_y(h)<deg_y(f_k)$
and $h\left(t^{\frac{v_0}{e_{k}}},\psi(t)\right)=bt^{\nu_k(h)}+\mbox{(h.o.t.)}$, we get $h(t^{v_0},\varphi(t))=bt^{\nu_f(h)}+\mbox{(h.o.t.)}$ with
$\nu_f(h)=e_{k}\nu_{k}(h)$.\end{proposition}
\Dem Considering $\varphi (t)=\displaystyle\sum_{i\geq v_0}a_it^i$, according with Corollary 5.3, Proposition 6.5 in \cite{popescu}, sections 6.2 and 8.1 in \cite{abramo}, for every $0< i<j<k\leq g$ we obtain:
$$\psi_k(t)-\phi\left(\zeta t^{\frac{e_{j}}{e_{k}}}\right)=\left \{\begin{array}{ll}
a_{\beta_{i}}\left(1-\zeta^{\frac{\beta_{i}}{e_{j}}}\right)t^{\frac{\beta_{i}}{e_{k}}}+ \mbox{(h.o.t.)} & \mbox{if}\ \zeta\in G_{i-1}\setminus G_{i};\\
a_{\beta_j}t^{\frac{\beta_{j}}{e_{k}}}+\mbox{(h.o.t.)} & \mbox{if}\ \zeta\in G_{j},
\end{array}
\right .$$
with $\zeta^{\frac{\beta_{i}}{e_{j}}}\neq 1$ if $\zeta\in G_{i-1}\setminus G_{i}$.

As $f_j(x,y)=\displaystyle\prod_{\overline{\zeta}\in\frac{G_0}{G_{j}}}\left ( y-\phi\left(\zeta x^{\frac{e_{j}}{v_0}}\right)\right )$ we have
$
f_j\left(t^{\frac{v_0}{e_{k}}},\psi(t)\right)=at^{\frac{v_j}{e_{k}}}+\mbox{(h.o.t.)}$,
with $a\neq 0$ and taking $k=g$ we get $f_j\left(t^{v_0},\varphi(t)\right)=at^{v_j}+\mbox{(h.o.t.)}$.

Moreover, given $h\in\mathbb{C}\{x\}[y]$ with $gr(h)<gr(f_k)$ by
Lemma \ref{decom-semirraizes} we can express $h=\sum_{
\alpha\in\mathbb{N}^{k}}b_{\alpha}f_0^{\alpha_0}\cdots
f_{k-1}^{\alpha_{k-1}}$ with $0\leq \alpha_j< n_{j+1}$ for all
$0\leq j< k$, $b_{\alpha}\in\mathbb{C}\{x\}$ and there exists a
unique $\gamma\in\mathbb{N}^{k}$ such that
$$\nu_f(h)=\nu_f(b_{\gamma}f_0^{\gamma_0}\cdots f_{k-1}^{\gamma_{k-1}})=\sum_{j=1}^{k-1}\gamma_{j-1}v_j+v_0\cdot ord_x(b_{\alpha})=e_k\nu_{f_k}(b_{\gamma}f_0^{\gamma_0}\cdots f_{k-1}^{\gamma_{k-1}}) =e_k\nu_{f_k}(h).$$
In addition, as
$f_j\left(t^{\frac{v_0}{e_{k}}},\psi(t)\right)=at^{\frac{v_j}{e_{k}}}+\mbox{(h.o.t.)}$
and $f_j\left(t^{v_0},\varphi(t)\right)=at^{v_j}+\mbox{(h.o.t.)}$
for all $0\leq j<k$, if
$h\left(t^{\frac{v_0}{e_{k}}},\psi(t)\right)=bt^{\nu_k(h)}+\mbox{(h.o.t.)}$,
then $h(t^{v_0},\varphi(t))=bt^{\nu_f(h)}+\mbox{(h.o.t.)}$.\cqd

Notice that $\Gamma_f=e_k\Gamma_k+\langle v_{k+1},\ldots ,v_g\rangle$ and
$$n_j=\frac{e_{j-1}}{e_j}=\frac{gcd(\beta_0,\ldots,\beta_{j-1})}{gcd(\beta_0,\ldots,\beta_j)}=\dfrac{gcd\left(\dfrac{\beta_0}{e_{k}},\ldots, \dfrac{\beta_{j-1}}{e_{k}} \right)}{gcd\left(\dfrac{\beta_0}{e_{k}},\ldots,\dfrac{\beta_j}{e_{k}} \right)}\ \ \ \mbox{for every}\ \ j=1,\ldots,k.$$
In this way, if $\mu_k$ denotes the Milnor number of $C_k$ we have $$\mu_f-1=\sum_{j=1}^{g}(n_j-1)v_j-v_0=e_{k}(\mu_k-1)+\sum_{j=k+1}^g(n_j-1)v_j.$$

The above remarks indicate how some topological data of a branch and
of their semiroots are related. More delicate question is: How we
can relate analytical invariants of $C_f$ and $C_k$? Investigating
this question is one of the aim of this work.

In the next section we will present results about some analytical
invariants of a plane branch associated with differential forms that
play an important role in the main results of this work.

\section{Differential forms on $C_f$}

Let $\Omega^1=\mathbb{C}\{x,y\}dx+\mathbb{C}\{x,y\}dy$ be the
$\mathbb{C}\{x,y\}$-module of holomorphic forms on $\mathbb{C}^2$.

Given $f\in\mathbb{C}\{x,y\}$ we denote $df=f_xdx+f_ydy$ and
$\mathcal{F}(f)=f\cdot\Omega^1+\mathbb{C}\{x,y\}\cdot df$. The
K\"ahler differential module of the branch $C_f$ is the
$\mathcal{O}_f$-module
$\Omega_f:=\Omega_{\mathcal{O}_f/\mathbb{C}}=\frac{\Omega^1}{\mathcal{F}(f)}$.
If $C_f$ is singular then $\Omega_f$ has a non trivial torsion
submodule $\mathcal{T}_f=\{\omega\in\Omega_f;\ h\omega=0\ \mbox{for
some}\ 0\neq h\in\mathcal{O}_f\}$ and
$\ell(\mathcal{T}_f)=\tau_f=\dim_{\mathbb{C}}\frac{\mathbb{C}\{x,y\}}{\langle
f,f_x,f_y\rangle}$ (see \cite{torsion} for instance).

Considering a parametrization $(x(t),y(t))\in\mathbb{C}\{t\}\times\mathbb{C}\{t\}$ of $C_f$ and the identification $\mathcal{O}_f\approx \mathbb{C}\{x(t),y(t)\}$ we define the $\mathcal{O}_f$-homomorphism
$$\begin{array}{cccc}
\Upsilon: & \Omega_f & \rightarrow & \mathbb{C}\{t\}\approx \overline{\mathcal{O}}_f \\
& A(x,y)dx+B(x,y)dy & \mapsto & t\cdot (A(x(t),y(t))\cdot
x'(t)+B(x(t),y(t))\cdot y'(t)),
\end{array}$$
and $\ker(\Upsilon)=\mathcal{T}_f$.

The set $\Lambda_f=\{\nu_f(\omega):=ord_t(\Upsilon(\omega));\
\overline{0}\neq\overline{\omega}\in\frac{\Omega_f}{\mathcal{T}_f}\}$
is an analytical invariant of $C_f$ and
$\Gamma_f\setminus\{0\}=\{\nu_f(dh);\ 0\neq
h\in\mathcal{O}_f\cap\langle t\rangle\}\subseteq \Lambda_f$. The set
$\Lambda_f$ has a prominent role in many works (see \cite{delorme}
and \cite{HefezHern} for instance) and it is related with the
Tjurina number $\tau_f$ of $C_f$ by the relation
$\tau_f=\mu_f-\sharp(\Lambda_f\setminus\Gamma_f)$. In \cite{algoritmo} we found an algorithm to compute the set $\Lambda_f$
for any irreducible curve plane or not and all possible set of values of differentials for plane curves with a fixed semigroup.

According to Saito (see \cite{Saito}) an element $W\in\frac{1}{f}\Omega^1$ is a logarithmic form along $C_f$ if there exist $\varpi\in\Omega^1$, $P,Q\in\mathbb{C}\{x,y\}, Q\not\in\langle f\rangle$ such that $QW=\frac{Pdf+f\varpi}{f}$. Saito denotes $\Omega(log\ C_f)=\{W;\ W\ \mbox{is a logarithmic form along}\ C_f\}$ and motivated by this, in the sequel, we will consider the $\mathbb{C}\{x,y\}$-module $f\cdot \Omega(log\ C_f)=\{f\cdot W;\ W\in\Omega(log\ C_f) \}\subset \Omega^1$. It follows that $\mathcal{T}_f=\frac{f\cdot\Omega(log\ C_f)}{\mathcal{F}(f)}$.

In section 2 of \cite{Saito}, Saito introduce the residue of a
logarithmic form that can be defined considering $f\cdot\Omega(log\
C_f)$. More precisely, given $\omega\in f\cdot\Omega(log\ C_f)$ with
$Q\omega=Pdf+f\varpi$ we call
$res(\omega)=\frac{\overline{P}}{\overline{Q}}\in
Frac(\mathcal{O}_f)$ the {\it residue of $\omega$} where
$Frac({\mathcal{O}_f})$ is the field of fraction of $\mathcal{O}_f$.
We have that $res(\omega)$ is well defined and
$\mathcal{R}_f:=\{res(\omega);\ \omega\in f\cdot\Omega(log\
C_f)\}\supseteq\overline{\mathcal{O}}_f$.

Pol (see Corollary 3.32, \cite{Pol}) shows that $\Lambda_f$
determines and it is determined by the set
$\Delta_f:=\{\nu_f(P)-\nu_f(Q);\
\frac{\overline{P}}{\overline{Q}}\in \mathcal{R}_f\}$. More
precisely, she proves that $\lambda\in\Lambda_f$ if and only if
$-\lambda\not\in\Delta_f$.

\begin{remark}
Given $\omega\in\Omega^1$ we denote $P_f(\omega):=\frac{\omega\wedge
df}{dx\wedge dy}\in J(f)$. In this way, we have that
$f\cdot\Omega(log\ C_f)=\{\omega\in\Omega^1;\ P_f(\omega)\in\langle
f\rangle\}$, that is, $\omega=Adx-Bdy\in f\cdot\Omega(log\ C_f)$ if
and only if $Af_y+Bf_x=P_f(\omega)=Mf$ for some
$M\in\mathbb{C}\{x,y\}$, or equivalently, $\frac{P_f(\omega)}{f}\in
(J(f):\langle f\rangle)$. In particular, we have
\begin{equation}\label{fomega}f_x\omega = Adf-fMdy\ \ \ \ \mbox{and}\ \ \ \ f_y\omega =-Bdf+fMdx,\end{equation}
    that is, $res(\omega)=\frac{A}{f_x}=\frac{B}{f_y}$.
    In addition, $\omega\in \mathcal{F}(f)$ if and only if $\frac{P_f(\omega)}{f}\in J(f)\subseteq (J(f):\langle f\rangle)$.
\end{remark}

The rest of this section is devoted to show that
$\nu_f((J(f):\langle f\rangle))$ and $\Delta_f$ (consequently,
$\Lambda_f$) are mutually determined and to revisit some results in
\cite{briancon} about equations $Af_y+Bf_x=C$.

As in the Section 2 we will assume that $f\in \mathbb{C}\{x\}[y]$ is an irreducible Weierstrass polynomial with $v_1=\nu_f(y)>\nu_f(x)=v_0=n=deg_y(f)=mult(f)>1$ and $v_0\nmid v_1$. In this way, we have $\nu_f(f_x)=\mu_f-1+v_1$, $\nu_f(f_y)=\mu_f-1+v_0=\sum_{i=1}^{g}(e_{i-1}-e_i)\beta_{i}$ (see Corollaries 7.15 and 7.16, \cite{abramo}) and
by (\ref{fomega})
$$\Delta_f=\{\nu_f(B)-(\mu_f-1+v_0);\ \omega=Adx+Bdy\in f\cdot\Omega(log\ C_f)\}.$$

\begin{remark}\label{jf-delta-partial} With the previous notation we have $\nu_f(J(f))=\left \{ \nu_f\left ( \frac{P_f(\omega)}{f}\right );\ \omega\in\mathcal{F}(f)\right \}$ and
$$\{\nu_f(res(\omega))\neq 0;\ \omega\in\mathcal{F}(f)\}=\Gamma_f\setminus\{0\}=\nu_f(J(f))-(\mu_f-1).
$$
In fact, as $\mathcal{F}(f)=f\cdot\Omega^1+\mathbb{C}\{x,y\}\cdot df$ it follows that
$\left \{ \frac{P_f(\omega)}{f};\ \omega\in\mathcal{F}(f)\right \}=J(f)$ and $\mathbb{C}\{x,y\}=\{res(\omega);\ \omega\in\mathcal{F}(f)\}$. We get the conclusion by (\ref{jac-gamma}).
\end{remark}

For the next result we introduce the following $\mathbb{C}$-vector spaces:
$$\begin{array}{ccc}
\vspace{0.2cm}
\mathcal{P}_s=\{h\in\mathbb{C}\{x\}[y];\ deg_y(h)<s\}; & &
\mathcal{E}(f)=\mathcal{P}_ndx+\mathcal{P}_{n-1}dy;\\
\mathcal{G}(f)=\mathbb{C}\{x,y\}df+\mathbb{C}\{x,y\}fdx; & & \mathcal{H}(f)=\mathbb{C}\{x\}df+f\Omega^1.
\end{array}$$

\begin{lemma}\label{direct} With the above notation we have
$\mathcal{E}(f)\oplus\mathcal{G}(f)=\Omega^1=\mathcal{E}(f)\oplus\mathcal{H}(f)$.
\end{lemma}
\Dem Given $\omega=Adx+Bdy\in\Omega^1$, by the Weierstrass division
theorem we can write $B=Qf_y+B_1$ with $B_1\in\mathcal{P}_{n-1}$ and
$A-Qf_x=Pf+A_1$ with $A_1\in\mathcal{P}_{n}$. So,
\begin{equation}\label{weierstrass}\omega=(A_1dx+B_1dy)+(Qdf+Pfdx)\in\mathcal{E}(f)+\mathcal{G}(f).\end{equation}
If $\omega\in\mathcal{E}(f)\cap\mathcal{G}(f)$ there exist
$A_1\in\mathcal{P}_n$, $B_1\in\mathcal{P}_{n-1}$,
$P,Q\in\mathbb{C}\{x,y\}$ such that $A_1dx+B_1dy=\omega=Qdf+Pfdx$,
that is, $A_1=Qf_x+Pf$ and $B_1=Qf_y$. As
$deg_y(f_y)=deg_y(f)-1=n-1$ it follows that $Q=P=0$, that is,
$\omega=0$. Hence, $\Omega^1=\mathcal{E}(f)\oplus\mathcal{G}(f).$

For the other equality we write $B=Q_1f+B_1$, $B_1=Q_0f_y+B_0$ and
$A-Q_0f_x=P_1f+A_0$ with $B_1, A_0\in\mathcal{P}_n$,
$B_0\in\mathcal{P}_{n-1}$ and
$Q_0\in\mathcal{P}_{1}=\mathbb{C}\{x\}$. So,
$$\omega=Adx+Bdy=(A_0dx+B_0dy)+(Q_0df+f(P_1dx+Q_1dy))\in\mathcal{E}(f)+\mathcal{H}(f)$$
and similarly to the previous case we conclude
$\mathcal{E}(f)\oplus\mathcal{H}(f)$. \cqd

In \cite{evelia} the authors call the expression (\ref{weierstrass})
is called the {\it Weierstrass form} of $\omega$ and it was
considered to study some properties of the foliation associated to
$\omega$.

Notice that if
$\omega=Qdf+Pfdx\in\mathcal{G}(f)\subset\mathcal{F}(f)\subset
f\cdot\Omega(log\ C_f)$ then $\nu_f(res(\omega))=\nu_f(Q)$ does not
have any relation with $\nu_f\left ( \frac{P_{f}(\omega)}{f}\right
)=\nu_f (Pf_y)$. On the other hand, we will show that for
$\omega\in\mathcal{E}(f)\cap f\cdot\Omega(log\ C_f)$ the values
$\nu_f(res(\omega))$ and $\nu_f\left ( \frac{P_{f}(\omega)}{f}\right
)$ are mutually determined.

In the sequel we follow the ideas of \cite{briancon}.

We denote the roots of $f\in\mathbb{C}\{x\}[y]$ by
$\{\varphi_i:=\varphi(\zeta^jx^{\frac{1}{n}});\ 1\leq j\leq
n\}\subset\mathbb{C}\{x^{\frac{1}{n}}\}$ where $\zeta$ is an $n$-th
primitive root of the unity, $\mathbb{K}:=Frac\left (
\mathbb{C}\{x^{\frac{1}{n}}\}\right )$ and
$\mathcal{P}_s^{\mathbb{K}}:=\{h\in\mathbb{K}[y];\ deg_y(h)<s\}$.

Let us consider the $\mathbb{K}$-basis $\mathbb{B}=\{\Phi_i:=\frac{f}{y-\varphi_i}=\prod_{j=1\atop j\neq i}^{n}(y-\varphi_j);\ i=1,\ldots ,n\}$ for $\mathcal{P}_n^{\mathbb{K}}$. Notice that $\Phi_i(\varphi_k)=\prod_{j=1\atop j\neq i}^{n}(\varphi_k-\varphi_j);\ i=1,\ldots ,n$, $$f_y(\varphi_k)=\Phi_k(\varphi_k)\ \ \ \mbox{and}\ \ \  f_x(\varphi_k)=-\varphi'_k\Phi_k(\varphi_k)=-\varphi'_kf_y(\varphi_k).$$ Moreover, given an element $h\in\mathcal{P}_n^{\mathbb{K}}$ we can express it using the basis $\mathbb{B}$ as $h=\sum_{i=1}^{n}h_i\Phi_i$ with $h_i=\frac{h(\varphi_i)}{\Phi_i(\varphi_i)}\in\mathbb{K}$.

Now, let us consider $\omega=Adx-Bdy\in \mathcal{E}(f)\cap
f\cdot\Omega(log\ C_f)$ with \begin{equation}\label{Mf}
P_f(\omega)=Bf_x+Af_y=Mf.
 \end{equation}
As $A\in\mathcal{P}_n$ and $B\in\mathcal{P}_{n-1}$ we get $M\in\mathcal{P}_{n-1}\subset\mathcal{P}_{n-1}^{\mathbb{K}}$.

The next result relates $M$ with $B$ by means the $\mathbb{K}$-basis $\mathbb{B}$.

\begin{lemma}[Lemma 2.2, \cite{briancon}] Given $M=\displaystyle\sum_{i=1}^{n}M_i\Phi_i$, $B=\displaystyle\sum_{i=1}^{n}B_i\Phi_i\in \mathcal{P}_{n-1}^{\mathbb{K}}$ and $A=\displaystyle\sum_{i=1}^nA_i\Phi_i\in\mathcal{P}_n^{\mathbb{K}}$ satisfying  (\ref{Mf}) we have $A_i=B_i\varphi_i^{\prime}$ and
\begin{equation}\label{Mi}
M_i=\left(\sum_{j=1\atop j\neq i}^n\left(\frac{\varphi_i^{\prime}-\varphi_j^{\prime}}{\varphi_i-\varphi_j} \right)\right)\cdot B_i-\sum_{j=1\atop j\neq i}^n\left(\left(\frac{\varphi_i^{\prime}-\varphi_j^{\prime}}{\varphi_i-\varphi_j}\right)\cdot B_j \right), \ \ \
\mbox{for every}\ \ \  i=1,\ldots,n.\end{equation}
\end{lemma}

As $\nu_f(M)=ord_t(M(t^n,\varphi_n(t)))$ and
$\nu_f(B)=ord_t(B(t^n,\varphi_i(t)))$ there exist
$b,c,d,u_j\in\mathbb{C}\setminus\{0\}$ such that
\begin{equation}\label{info}\begin{array}{lll}
B(x,\varphi_i)=b.\zeta^{i\nu_f(B)}x^{\frac{\nu_f(B)}{n}}+\mbox{(h.o.t.)},
& &
\Phi_i(\varphi_i)=f_y(x,\varphi_i)=c.\zeta^{i\nu_f(f_y)}x^{\frac{\nu_f(f_y)}{n}}+\mbox{(h.o.t)},\\
M(x,\varphi_n)=d.x^{\frac{\nu_f(M)}{n}}+\mbox{(h.o.t.)},& &
\varphi_n-\varphi_j=u_jx^{\frac{m_j}{n}}+\mbox{(h.o.t.)}\ \ (j\neq
n),\end{array} \end{equation} where $\zeta$ is an $n$-th primitive
root of the unity and $m_j=\beta_k$ if $\zeta^j\in G_{k-1}\setminus
G_k$. So, we get
$-\left(\frac{\varphi_n^{\prime}-\varphi_j^{\prime}}{\varphi_n-\varphi_j}\right)=-\frac{m_j}{nx}+\mbox{(h.o.t.)}$.

Considering $i=n$ in (\ref{Mi}) and expanding the expressions we obtain
\begin{equation}\label{eqmb}d\cdot x^{\frac{\nu_f(M)-\nu_f(f_y)}{n}}+(\mbox{h.o.t.})=
\sum_{i=1}^{n-1}m_i\left ( 1 -\zeta^{i(\nu_f(B)-\nu_f(f_y))}\right
)\cdot \dfrac{b}{n}\cdot
x^{\frac{\nu_f(B)-\nu_f(f_y)-n}{n}}+(\mbox{h.o.t.}). \end{equation}

Denote $S:=\sum_{i=1}^{n-1}m_i\left ( 1
-\zeta^{i(\nu_f(B)-\nu_f(f_y))}\right )$. As, $m_i=\beta_k$ for
$\eta=\zeta^{i}\in G_{k-1}\setminus G_{k}$ we get
\begin{equation}\label{S}
S=\sum_{i=1}^{g}\left (\sum_{\eta\in G_{i-1}\setminus
G_{i}}\left(1-\eta^{\nu_f(B)-\nu_f(f_y)}\right
)\right)\beta_{i}.\end{equation}

The next two lemmas allow us to conclude that $S\neq 0$, consequently we obtain $\nu_f(M)=\nu_f(B)-n$ and $d=\frac{S\cdot b}{n}$.

\begin{lemma}\label{kexiste} If $B\in\mathcal{P}_{n-1}$ with $\nu_f(B)=\displaystyle\sum_{i=0}^gs_iv_i$, $0\leq s_i<n_i$ and $s_0\geq 0$ then there exists $k\in G_0$ such that $e_k\nmid\alpha:=\displaystyle\sum_{i=1}^g(n_i-1-s_i)v_i$.\end{lemma}
\Dem Let us consider a complete system of semiroots
$\{f_0,\ldots,f_{g} \}$ of $f$. As $deg_y(B)<n-1$, according Lemma
\ref{decom-semirraizes} we can write
$$B=\displaystyle\sum_{finite\atop r\in\mathbb{N}^{g+1}}a_{r}f_0^{r_1}\cdots f_{g-1}^{r_{g}}\ \ \mbox{with}\ \
\nu_f(B)=\nu_f(a_{s}f_0^{s_1}\cdots f_{g-1}^{s_g}
)=\sum_{i=0}^{g}s_iv_i.$$

Let us suppose that $e_k\mid \alpha$ for every $k\in
G_0\setminus\{1\}$. If there exists $k=\max\{i>0;\ s_i< n_i-1\}$
then, by (\ref{vi}) and (\ref{milnor}), we have
$(n_i-1-s_i)v_i>n_{i-1}v_{i-1}>e_{i-1}(\mu_{i-1}-1)$. But this
implies $(n_i-1-s_i)v_{i}\in\langle v_0,\ldots ,v_{i-1}\rangle$ that
is a contradiction because $n_i=\min\{j>0;\ j.v_i\in\langle
v_0,\ldots ,v_{i-1}\rangle\}$. So, if $e_k\mid \alpha$ for every
$k\in G_0\setminus\{1\}$ then $s_k=n_k-1$ for every $1\leq k\leq g$.
By Remark \ref{degree} we would have $deg_y(B)= n-1$ that is an
absurd. \cqd

\begin{lemma}\label{casos} Given $\alpha\in\mathbb{Z}$ and $k=\max\{i\in G_0;\ e_i\nmid\alpha\}$ we have
$$\sum_{\eta\in G_{j}\setminus G_{j+1}}\eta^{\alpha}=\left\{\begin{array}{lll}e_{j}-e_{j+1}&\mbox{if}& k+1\leq j< g\\
-e_{j+1}&\mbox{if}&j=k\\
0&\mbox{if}&0\leq j<k.\end{array} \right. $$
\end{lemma}
\Dem As $\sharp G_k=e_k$  we have that $\sum_{\eta\in
G_k}\eta^{\alpha}=\left\{\begin{array}{lll}
0&\mbox{if}&e_k\nmid \alpha\\
e_k&\mbox{if}&e_k\mid \alpha.
\end{array}\right.$

Therefore, the result follows considering $\sum_{\eta\in
G_{j}}\eta^{\alpha}=\sum_{\eta\in G_{j}\setminus
G_{j+1}}\eta^{\alpha}+\sum_{\eta\in G_{j+1}}\eta^{\alpha}$ and
recalling that $e_{j}\mid e_{k}$ for $j\leq k$. \cqd

As a consequence of the previous two lemmas we obtain the following result.

\begin{proposition}\label{BM1} If $\omega=Adx-Bdy\in\mathcal{E}(f)\cap f\cdot\Omega(log\ C_f)$  with $P_f(\omega)=Bf_x+Af_y=Mf$ then
$$ M(t^n,\varphi_n(t))=\dfrac{e_kv_{k+1}b}{n}t^{\nu_f(M)}+(\mbox{h.o.t.}),\ \ \ \
B(t^n,\varphi_n(t))=bt^{\nu_f(B)}+(\mbox{h.o.t.})$$
and $\nu_f(B)=\nu_f(M)+v_0$, where $k=\max\{i\in G_0;\ e_i\nmid \nu_f(B)-\nu_f(f_y)\}$.
\end{proposition}
\Dem Remark that Lemma \ref{kexiste} give us the existence of
$k=\max\{i\in G_0;\ e_i\nmid \nu_f(B)-\nu_f(f_y)\}$. By (\ref{eqmb})
it is sufficient to show that in (\ref{S}) we have
$S=\sum_{i=1}^{g}\left (\sum_{\eta\in G_{i-1}\setminus
G_{i}}\left(1-\eta^{\nu_f(B)-\nu_f(f_y)}\right )\right)\beta_{i}\neq
0$. But this is a consequence of Lemma \ref{casos} and (\ref{vi})
because
$$ S=\sum_{i=1}^{g}(e_{i-1}-e_i)\beta_i-\left(-e_{k+1}\beta_{k+1}+\sum_{j=k+2}^{g}(e_{j-1}-e_{j})\beta_{j}  \right)=\sum_{i=1}^{k}(e_{i-1}-e_i)\beta_i+e_k\beta_{k+1}=e_kv_{k+1}\neq 0.$$
\cqd

The next result shows that $\Delta_f$ and $\nu_f(J(f):f)$ are
mutually determined.

\begin{theorem}\label{conjuntos} Let $f\in\mathbb{C}\{x\}[y]$ be an irreducible Weierstrass polynomial with $deg_y(f)=v_0=n$ and $\Gamma_f=\langle v_0,v_1,\ldots ,v_g\rangle$, then $\{\delta\neq 0;\ -\delta\not\in\Lambda_f\}=\Delta_f\setminus\{0\}=\nu_f(J(f):f)-(\mu_f-1)$.
\end{theorem}
\Dem By Pol (see Corollary 3.32, \cite{Pol}) it follows that $\{\delta;\ -\delta\not\in\Lambda_f\}=\Delta_f$ and by Remark (\ref{jf-delta-partial}) we have that $\Gamma_f\setminus\{0\}\subseteq \Delta_f\cap (\nu_f(J(f):f)-(\mu_f-1))$.

Given $\alpha\in\Delta_f\setminus\Gamma_f$ there exists
$\omega=Adx+Bdy\in f\cdot\Omega(log\ C_f)$ with
$\nu_f(res(\omega))=\nu_f\left ( \frac{B}{f_y}\right )=\alpha$ and,
according to Lemma \ref{direct}, we can write
\begin{equation}\label{exp}
\omega=Adx+Bdy=A_0dx+B_0dy+Rdf+(Pdx+Qdy)f\in\mathcal{E}(f)\oplus\mathcal{H}(f),
\end{equation} with $A_0dx+B_0dy\in\mathcal{E}(f)\cap f\cdot\Omega(log\ C_f)$
and $R\in\mathbb{C}\{x\}$. In particular,
$\alpha+\mu_f-1+v_0=\nu_f(B)=\nu_f(B_0+Rf_y+Qf)=\nu_f(B_0+Rf_y)$ and
$\nu_f(R)=s_0v_0$.

If $\nu_f(Rf_y)<\nu_f(B_0)$ we should have $\alpha+\mu_f-1+v_0=\nu_f(Rf_y)=s_0v_0+\mu_f-1+v_0$, i.e., $\alpha=s_0v_0\in\Gamma_f$ that is an absurd.

If $\nu_f(B_0)=\nu_f(Rf_y)$, then $\nu_f(B_0)=s_0v_0+\sum_{i=1}^g(n_i-1)v_i$. But, in this way, Lemma \ref{degree} gives us $deg_y(B_0)\geq n-1$ and we get a contradiction.

In this way $\nu_f(Rf_y)>\nu_f(B_0)=\nu_f(B)=\alpha+\mu_f-1+v_0$ and taking $\omega_0=A_0dx+B_0dy\in\mathcal{E}(f)\cap f\cdot\Omega (log\ C_f)$ with $P_f(\omega_f)=M_0f$ we have $M_0\in (J(f):f)$, by Proposition \ref{BM1},
$\alpha+\mu_f-1=\nu_f(B)-v_0=\nu_f(B_0)-v_0=\nu_f(M_0)\in\nu_f(J(f):f)$.

Reciprocally, consider $\alpha+\mu_f-1=\nu_f(M)\in\nu_f(J(f):f)$ with $\alpha\not\in\Gamma_f$.

According to Saito, we have $\overline{\mathcal{O}}_f\subseteq\mathcal{R}_f$ so, if $\alpha>0$ then $\alpha\in\Delta_f$.

On the other hand, if $\alpha<0$ then, by  Lemma \ref{direct}, there
exists \begin{equation}\label{exp2}
\omega=Adx+Bdy=A_0dx+B_0dy+Qdf+Pfdx\in f\cdot\Omega (log\ C_f),
\end{equation} with $\omega_0=A_0dx+B_0dy\in\mathcal{E}(f)\cap f\cdot\Omega
(log\ C_f)$, $P,Q\in\mathbb{C}\{x,y\}$ and
$Mf=P_f(\omega)=P_f(\omega_0)+Pf_yf=(M_0+Pf_y)f$. As $\alpha<0$ we
have $\alpha+\mu_f-1+v_0<\nu_f(f_y)\leq \nu_f(Pf_y)$ and, by
Proposition \ref{BM1}, $\nu_f(M)=\nu_f(M_0)=\nu_f(B_0)-v_0$, that
is, $\nu_f(B_0)=\alpha+\mu_f-1+v_0<\nu_f(f_y)$. In this way,
$\nu_f(B)=\nu_f(B_0+Qf_y)=\nu_f(B_0)$ and
$\nu_f(res(\omega))=\nu_f\left(\dfrac{B}{f_y}
\right)=\nu_f\left(\dfrac{B_0}{f_y} \right)=\alpha\in\Delta_f$. \cqd

As a consequence of the considerations in the proof of the above theorem we have the following result:

\begin{corollary}\label{ordem-delta}
    For every $\alpha\in\Delta_f\setminus\Gamma_f$ there exists $\omega\in\mathcal{E}(f)\cap f\cdot\Omega (log\ C_f)$ such that $\nu_f(res(\omega))=\alpha=\nu_f\left ( \frac{P_f(\omega)}{f}\right )-(\mu_f-1)$.
\end{corollary}

Proposition \ref{BM1} relates $\nu_f(B)$ (equivalently $\nu_f(A)$)
and $\nu_f(M)$ for $A\in\mathcal{P}_n$ and $B,
M\in\mathcal{P}_{n-1}$ satisfying $Af_y+Bf_x=Mf$. As $f_x$ and $f_y$
are coprime, given $H\in \mathbb{K}[y]$ there exist
$A,B\in\mathbb{K}[y]$ such that $Af_y+Bf_x=H$. In \cite{briancon},
Brian\c con, Maisonobe and Torrelli consider such equation and, with
the notations introduced before, obtain the following result:

\begin{lemma}[Proposition 2.16, \cite{briancon}] If $B=\displaystyle\sum_{i=1}^nB_i\Phi_i,H=\displaystyle\sum_{i=1}^nH_i\Phi_i\in\mathcal{P}_{n-1}^K$ and $A=\displaystyle\sum_{i=1}^nA_i\Phi_i\in\mathcal{P}_n^K$ satisfy $Af_y+Bf_x=H$ then $A=\sum_{i=1}^n\left(B_i\varphi^{\prime}_i+\dfrac{H_i}{\Phi_i(\varphi_i)} \right)\Phi_i$ and
\begin{eqnarray*}\left(\sum_{j=1\atop j\neq i}^n\dfrac{\varphi^{\prime}_i
-\varphi^{\prime}_j}{\varphi_i-\varphi_j} \right)\cdot B_i-\sum_{j=1\atop j\neq i}^n\left(\left(\dfrac{\varphi^{\prime}_i-\varphi^{\prime}_j}{\varphi_i-\varphi_j}\right)\cdot B_j \right)&=&\left(\sum_{j=1\atop j\neq i}^n\dfrac{\Phi_j(\varphi_j)-\Phi_i(\varphi_i)}{(\varphi_i-\varphi_j)\Phi_i(\varphi_i)\Phi_j(\varphi_j)} \right)\cdot H_i\\&-&\sum_{j=1\atop j\neq i}^n\left(\left(\dfrac{\Phi_j(\varphi_j)-\Phi_i(\varphi_i)}{(\varphi_i-\varphi_j)\Phi_i(\varphi_i)\Phi_j(\varphi_j)}\right)\cdot H_j \right). \end{eqnarray*}
\end{lemma}

Using the computations presented in Proposition \ref{BM1} we can
relate $\nu_f(H)$ and $\nu_f(B)$.

\begin{proposition}\label{BM2} With the previous notation we get
$$\nu_f(H)\leq \nu_f(B)+\mu_f-1+\beta_g$$
where $\beta_g$ is the greatest characteristic exponent of
$C_f$.\end{proposition} \Dem By (\ref{Mi}), (\ref{info}),
(\ref{eqmb}) and (\ref{S}) we have
\begin{equation}\label{EQ1}\left(\sum_{j=1}^{n-1}\dfrac{\varphi^{\prime}_n-\varphi^{\prime}_j}{\varphi_n-\varphi_j}
\right)\cdot B_n-\sum_{j=1}^{n-1}\left(\left(\dfrac{\varphi^{\prime}_n-\varphi^{\prime}_j}{\varphi_n-\varphi_j}\right)\cdot
B_j \right)=
\frac{bS}{cn}x^{\frac{\nu_f(B)-\nu_f(f_y)-n}{n}}+\mbox{(h.o.t.)};\end{equation}
$$\dfrac{\Phi_n(\varphi_n)-\Phi_j(\varphi_j)}{(\varphi_n-\varphi_j)\Phi_n(\varphi_n)\Phi_j(\varphi_j)}
=\left(\dfrac{1-\eta^{j\nu_f(f_y)}}{cu_j\eta^{j\nu_f(f_y)}}\right)x^{\frac{-(\nu_f(f_y)+m_j)}{n}}+\mbox{(h.o.t.)},
$$
where $\eta$ is a $n$-th primitive root of the unity,
$m_j\in\{\beta_1,\ldots ,\beta_g\}$ and $S\neq 0$ as showed in
Proposition \ref{BM1}.

Denoting
$H(x,\varphi_j)=h\eta^{j\nu_f(H)}x^{\frac{\nu_f(H)}{n}}+\mbox{(h.o.t.)}$
with $h\in\mathbb{C}\setminus\{0\}$ we obtain {\small
$$\hspace{-3cm}\left(\sum_{j=1}^{n-1}\dfrac{\Phi_j(\varphi_j)-\Phi_n(\varphi_n)}{(\varphi_n-\varphi_j)\Phi_n(\varphi_n)\Phi_j(\varphi_j)}
\right)\cdot H_n-\sum_{j=1}^{n-1}\left(\left(\dfrac{\Phi_j(\varphi_j)-\Phi_n(\varphi_n)}{(\varphi_n-\varphi_j)\Phi_n(\varphi_n)\Phi_j(\varphi_j)}\right)\cdot
H_j \right )=$$
\begin{equation}\label{EQ2}\hspace{5cm}=\sum_{j=1}^n\left(\frac{1-\eta^{j\nu_f(f_y)}}{u_j\eta^{2j\nu_f(f_y)}}\right)\left(\dfrac{h}{c^2}\eta^{j\nu_f(H)}\right)x^{\frac{\nu_f(H)-2\nu_f(f_y)-m_j}{n}}+\mbox{(h.o.t.)}.
\end{equation} }

By the previous lemma, expressions (\ref{EQ1}) and (\ref{EQ2}) there exists
$1\leq j\leq n$ such that
$$\dfrac{\nu_f(H)-2\nu_f(f_y)-m_j}{n}\leq \dfrac{\nu_f(B)-\nu_f(f_y)-n}{n},$$
that is, $\nu_f(H)\leq \nu_f(B)+\mu_f-1+m_j\leq
\nu_f(B)+\mu_f-1+\beta_g$. \cqd

\section{Relating analytical invariants to $C_k$ and $C_f$}

Given a $k$-semiroot $f_k$ of $f$ with $0\leq k<g$, by Lemma
\ref{direct}, we get
$\Omega_1=\mathcal{E}(f_k)\oplus\mathcal{G}(f_k)$ with
$\mathcal{G}(f_k)\subseteq\mathcal{F}(f_k)$. In this way, for each
$\varpi\in\Omega_{f_k}=\frac{\Omega^1}{\mathcal{F}(f_k)}$ there
exist $\omega\in\mathcal{E}(f_k)$ such that
$\overline{\omega}=\varpi$. In this section we will analyze
$\nu_f(\omega)$ for $\omega\in\mathcal{E}(f_k)$ and we will describe
elements in $\Lambda_f\setminus\Gamma_f$ from elements of
$\Lambda_k\setminus\Gamma_k$.

Notice that for $k=0$, that is $f_0=y+a(x)\in\mathbb{C}\{x\}[y]$, we
have $\mathcal{E}(f_0)=\mathbb{C}\{x\}dx+\mathbb{C}\{x\}dy$ and
$\nu_f(\mathcal{E}(f_0))\subseteq\{v_0+kv_0,v_1+kv_0;\
k\in\mathbb{N}\}\subset\Gamma_f$.

In what follows we consider $1\leq k<g$.

Given $\omega=Adx-Bdy\in\mathcal{E}(f_k)$, that is, $A, B\in\mathbb{C}\{x\}[y]$ with $deg_y(A)<deg_y(f_k)-1=\frac{n}{e_k}-1$ and $deg_y(B)<\frac{n}{e_k}$ we can write
\begin{equation}\label{basico}
P_{f_k}(w)=A\cdot (f_k)_y+B\cdot (f_k)_x=H+Mf_k
\end{equation} with $H, M\in\mathbb{C}\{x\}[y]$, $deg_y(H)<\frac{n}{e_k}$ and $deg_y(M)<\frac{n}{e_k}-1$.

We will consider the cases: $\omega\in\mathcal{E}(f_k)\cap
f_k\cdot\Omega (log\ C_k)$ or $\omega\in\mathcal{E}(f_k)\setminus
f_k\cdot\Omega (log\ C_k)$.

\subsection{Case: $\omega\in\mathcal{E}(f_k)\cap f_k\cdot\Omega (log\ C_k)$}

In this case, $H=0$ in (\ref{basico}) and we obtain the following
result:

\begin{proposition}\label{torcao} If $\omega=Adx-Bdy\in\mathcal{E}(f_k)\cap f_k\cdot
\Omega(log\ C_k)$ with $P_{f_k}(w)=Mf_k$ then
$$\nu_f(\omega)=v_{k+1}-e_k(\mu_k-1-\nu_{k}(M))=v_{k+1}-e_k\left(\mu_k-1+\frac{v_0}{e_k}-\nu_{f_k}(B)\right).$$
In particular, $\nu_f\left(\mathcal{E}(f_k)\cap f_k\cdot\Omega(log\
C_k) \right)\subseteq\{v_{k+1}-e_k\delta>0; \
0\neq\delta\in\mathbb{Z}\setminus\Lambda_k\}.$ \end{proposition}
\Dem If $\left( t^{\frac{v_0}{e_k}},\psi(t)\right )$ and
$(t^{v_0},\varphi(t))$ denote parametrizations for $C_{k}$ and $C_f$
respectively, then by Proposition \ref{aval-semiraiz} and
Proposition \ref{BM1}, there exists $0\leq j\leq k-1$ such that
$$\begin{array}{ll} M\left(t^{\frac{v_0}{e_k}},\psi(t)\right)=\frac{e_{j}v_{j+1}}{n}bt^{\nu_{k}(M)}+\mbox{(h.o.t.)}; & M\left(t^{v_0},\varphi(t)\right)=\frac{e_{j}v_{j+1}}{n}bt^{\nu_{f}(M)}+\mbox{(h.o.t.)}; \\
B\left(t^{\frac{v_0}{e_k}},\psi(t)\right)=bt^{\nu_k(B)}+\mbox{(h.o.t.)}; &
B\left(t^{v_0},\varphi(t)\right)=bt^{\nu_f(B)}+\mbox{(h.o.t.)}
\end{array}$$
with $\nu_f(B)=e_k\nu_k(B)=e_k\left (\nu_k(M)+\frac{v_0}{e_k}\right )=\nu_f(M)+v_0$.

As $deg_y((f_k)_y)<deg_y(f_k)$, Proposition \ref{aval-semiraiz} gives us  $\nu_f((f_k)_y)=e_k\nu_k((f_k)_y)=e_k(\mu_k-1)+v_0$.

By (\ref{fomega}) we get $(f_k)_y\omega=Mf_kdx-Bdf_k$. If
$f_k(t^{v_0},\varphi(t))=at^{v_{k+1}}+\mbox{(h.o.t.)}$ with
$a\in\mathbb{C}\setminus\{0\}$ then
$(Mf_kdx-Bdf_k)(t^{v_0},\varphi(t))=ab(e_jv_{j+1}-v_k)t^{\nu_f(B)}+\mbox{(h.o.t.)}$
and
$$v_f(\omega)=
v_{k+1}-e_k\left(\mu_k-1-\nu_{k}(M)\right)=v_{k+1}-e_k\left(\mu_{k}-1+\dfrac{v_0}{e_k}-\nu_{f_k}(B)\right).$$

In particular, by Theorem \ref{conjuntos},
$-\delta=\mu_k-1-\nu_k(M)\in\Delta_k\setminus\{0\}$, or
equivalently, $\delta\not\in\Lambda_k$. So,
$\nu_f\left(\mathcal{E}(f_k)\cap f_k\cdot\Omega(log\ C_k)
\right)\subseteq\{v_{k+1}-e_k\delta>0; \
0\neq\delta\in\mathbb{Z}\setminus\Lambda_k\}$. \cqd

As a consequence we obtain the following corollary:

\begin{corollary}\label{valores1} Considering $\Gamma_f=\langle v_0,\ldots ,v_g\rangle$
and maintaining the above notation we have
\begin{equation}\label{log-value}\nu_f\left(\mathcal{E}(f_k)\cap f_k\cdot\Omega(log\ C_k)
\right)\ \cap \ (\Lambda_f\setminus\Gamma_f)=\{v_{k+1}-e_{k}\delta;
\ \delta\in\mathbb{N}^{\ast}\setminus\Lambda_k \ \mbox{or} \
-\delta\in\mathbb{N}\setminus\Gamma_k \}. \end{equation} Moreover,
for every $\delta\in \{\alpha; \
\alpha\in\mathbb{N}^{\ast}\setminus\Lambda_{k}\}\dot\cup\{-\alpha; \
\alpha\in\mathbb{N}\setminus\Gamma_k \}$ we have
$$\sum_{i=k+1}^gs_iv_i-e_{k}\delta\in\Lambda_f\setminus\Gamma_f $$
where $0\leq s_j<n_{j}$, $k< j\leq g$ with $s_{k+1}\neq 0$.
\end{corollary}
\Dem Firstly, remark that for $-\delta\in\Gamma_k$ we have
$v_{k+1}-e_k\delta\in\Gamma_f$, then by the above proposition we get
$\nu_f\left(\mathcal{E}(f_k)\cap f_k\cdot\Omega(log\ C_k) \right)\
\cap \ (\Lambda_f\setminus\Gamma_f)\subseteq\{v_k-e_{k}\delta; \
\delta\in\mathbb{N}^{\ast}\setminus\Lambda_k \ \mbox{or} \
-\delta\in\mathbb{N}\setminus\Gamma_k \}$.

On the other hand, if $\delta\in\mathbb{N}^*\setminus\Lambda_{k}$,
respectively  $-\delta\in\mathbb{N}\setminus\Gamma_k$ then, by
Theorem \ref{conjuntos} (applied for $f_k$), there exists
$\omega_0=A_0dx+B_0dy\in \mathcal{E}(f_k)\cap f_k\cdot\Omega(log\
C_k)$ (see (\ref{exp}), respectively (\ref{exp2})) such that
$P_{f_k}(\omega_0)=M_0f_k$ with $\nu_k(M_0)=\mu_k-1-\delta$,
consequently by the previous proposition we have
$v_{k+1}-e_k\delta=v_{k+1}-e_k(\mu_k-1-\nu_k(M_0))=\nu_k(\omega_0)\in\Lambda_f\cap\nu_f\left(\mathcal{E}(f_k)\cap
f_k\cdot\Omega(log\ C_k) \right)$.

Moreover, if $\{f_0,\ldots ,f_g\}$ is a complete system of semiroots
of $f$, then $$\gamma=\sum_{i=k+1}^gs_iv_i-e_{k}\delta=\nu_f\left (
f_{k}^{s_{k+1}-1}\prod_{i=k+1}^{g-1}f_i^{s_{i+1}}\omega_0\right )\in\Lambda_f$$ for $0\leq s_{i}<n_{i}$, $k< i\leq g$ with
$s_{k+1}\neq 0$. So, to obtain the proposition it is sufficient to show that
$\gamma\not\in\Gamma_f$ for every  $\delta\in \{\alpha; \
\alpha\in\mathbb{N}^{\ast}\setminus\Lambda_{k}\}\dot\cup\{-\alpha; \
\alpha\in\mathbb{N}\setminus\Gamma_k \}$. Notice that, taking $s_{k+1}=1,s_{k+2}=\ldots =s_{g-1}=0$ we get (\ref{log-value}).

If $-\delta\in\mathbb{N}\setminus\Gamma_k$ then, by Remark \ref{padrao},
we can write
$-\delta=\sum_{i=0}^{k}s_i\frac{v_i}{e_k}$ with
$0\leq s_i<n_i$; $1\leq i\leq k$ and $s_0<0$, consequently
$\gamma=\sum_{i=0}^{g}s_iv_i\not\in\Gamma_f$.

If $\delta\in\mathbb{N}^*\setminus\Lambda_k$ and $v_{k+1}-e_k\delta\in\Gamma_f$, as $v_{k+1}-e_k\delta<v_{k+1}$ we must have that $e_k$ divides $v_{k+1}-e_k\delta$, that is, $e_k\mid v_{k+1}$ that is an absurd. So, by Remark \ref{padrao}, $v_{k+1}-e_k\delta=\sum_{i=0}^{k}s_iv_i$ with $0\leq s_i<n_i$; $1\leq i\leq k$ and $s_0<0$ and
$\gamma=\sum_{i=0}^{k}s_iv_i+(s_{k+1}-1)v_{k+1}+\sum_{i=k+2}^{g}s_iv_i\not\in\Gamma_f$.
\cqd

\subsection{Case: $\omega\in\mathcal{E}(f_k)\setminus f_k\cdot\Omega (log\ C_k)$}

With the same above notations we consider
$\omega=Adx-Bdy\in\mathcal{E}(f_k)\setminus f_k\cdot\Omega (log\
C_k)$  and $\nu_k(\omega)\in\Lambda_k$. Remark that
$P_{f_k}(\omega)=A\cdot (f_k)_y+B\cdot (f_k)_x=H+Mf_k$ and
\begin{equation}\label{19}
(f_k)_y\omega=P_{f_k}(\omega)dx-Bdf_k=Hdx+Mf_kdx-Bdf_k,
\end{equation}
with $B, H, M\in\mathbb{C}\{x\}[y]$,
$deg_y(H)<deg_y(f_k)=\frac{n}{e_k}$ and $deg_y(B),
deg_y(M)<\frac{n}{e_k}-1$. As $\nu_k((f_k)_y)=\mu_k-1+\frac{n}{e_k}$
we get
\begin{equation}\label{valor}\nu_k(\omega)=\nu_k(H)-(\mu_k-1).\end{equation}

On the other hand, as $(f_k)_y$ and $(f_k)_x$ are coprime in
$\mathbb{C}((x))[y]$ there exist $A',B'\in\mathbb{C}((x))[y]$ with
$deg_y(A')<\frac{n}{e_k}$ and $deg_y(B')<\frac{n}{e_k}-1$ such that
\begin{equation}\label{linha}
A'\cdot (f_k)_y+B'\cdot (f_k)_x=Mf_k. \end{equation} In particular, there
exists $Q\in\mathbb{C}\{x\}$ such that $A_0=QA',
B_0=QB'\in\mathbb{C}\{x\}[y]$ satisfy
$A_0\cdot (f_k)_y+B_0\cdot (f_k)_x=QMf_k$ and, by Proposition \ref{BM1},
$\nu_k(B_0)=\nu_k(QM)+\frac{n}{e_k}$ with
$B_0(t^{\frac{v_0}{e_k}},\varphi_k(t))=bt^{\nu_k(B_0)}+
\mbox{(h.o.t.)}$,
$QM(t^{\frac{v_0}{e_k}},\varphi_k(t))=\frac{e_jv_jb}{v_0}t^{\nu_k(QM)}+
\mbox{(h.o.t)}$ where $e_j=\min\{e_i;\ e_i\nmid
\nu_k(B_0)-\nu_k((f_k)_y)\}<e_k$ and
$(t^{\frac{v_0}{e_k}},\varphi_k(t))$ is a parametrization of $f_k$.
Consequently, by Proposition \ref{aval-semiraiz}, we have
\begin{equation}\label{21}\nu_f(B')+v_{k+1}=\nu_f(M)+v_0+v_{k+1}=\nu_f(Mf_kdx-B'df_k).
\end{equation}

In this way, $(A-A')\cdot (f_k)_y+(B-B')\cdot (f_k)_x=H$ and, by Proposition
\ref{BM2} $\nu_k(H)\leq \nu_k(B-B')+\mu_k-1+\frac{\beta_k}{e_k}$, or
equivalently, by Proposition \ref{aval-semiraiz}, $\nu_f(H)+v_0\leq
\nu_f(B-B')+e_k(\mu_k-1)+\beta_k<\nu_f(B-B')+e_k(\mu_k-1)+\beta_{k+1}.$
Using (\ref{milnor}) and (\ref{vi}) for $f_k$ we obtain
$$\nu_f(H)+v_0<\nu_f(B-B')+v_{k+1}.$$
Notice that $\nu_f(H)+v_0=e_k\left ( \nu_k(H)+\frac{v_0}{e_k}\right
)$ and $e_{k}\nmid (\nu_f(B)+v_{k+1})$ so, $\nu_f(H)+v_0\neq
\nu_f(B)+v_{k+1}$.

\begin{lemma}\label{4.3}
Considering (\ref{19}) and the above notations we have that
$$\nu_f(\omega)=\left \{ \begin{array}{ll}
e_k\nu_k(\omega) & \mbox{if}\ \ \ \nu_f(H)+v_0<\nu_f(B)+v_{k+1};\\
e_k\nu_k(B)+\beta_{k+1} & \mbox{if}\ \ \
\nu_f(H)+v_0>\nu_f(B)+v_{k+1}.
\end{array}\right .$$
\end{lemma}
\Dem By (\ref{19}) we have that
$$\nu_f(\omega)=\nu_f(Hdx+Mf_kdx-Bdf_k)-\nu_f((f_k)_y)=\nu_f(Hdx+Mf_kdx-Bdf_k)-e_k\left (\mu_k-1+\frac{v_0}{e_k}\right).$$

Suppose $\nu_f(H)+v_0<\nu_f(B)+v_{k+1}$. If
$\nu_f(B-B')=\min\{\nu_f(B), \nu_f(B')\}$, then
$\nu_f(B-B')+v_{k+1}\leq \nu_f(B')+v_{k+1}=\nu_f(M)+v_0+v_{k+1}$. If
$\nu_f(B-B')>\min\{\nu_f(B), \nu_f(B')\}$, then by (\ref{21}) we get
$\nu_f(Mf_kdx-Bdf_k)=\nu_f(B)+v_{k+1}$ and by (\ref{valor})
$\nu_f((f_k)_y\omega)=\nu_f(Hdx)=e_k\nu_k((f_k)_y\omega)$.

Now suppose that $\nu_f(H)+v_0>\nu_f(B)+v_{k+1}$. As
$\nu_f(H)+v_0<\nu_f(B-B')+v_{k+1}$ we must have
$\nu_f(B-B')>\min\{\nu_f(B), \nu_f(B')\}$, then by (\ref{21}) we get
$\nu_f(Mf_kdx-Bdf_k)=\nu_f(B)+v_{k+1}$ and
$\nu_f(\omega)=\nu_f(Bdf_k)-e_k\left ( \mu_k-1+\frac{v_0}{e_k}\right
)=v_{k+1}-e_k\left ( \mu_k-1+\frac{v_0}{e_k}-\nu_k(B)\right
)=\beta_{k+1}+e_k\nu_k(B)$ where the last equality follows because
$e_k(\mu_k-1)+v_0=n_kv_k-\beta_k=v_{k+1}-\beta_{k+1}$. \cqd

\begin{remark}\label{remark}
Notice that by the previous result if
$\nu_f(H)+v_0<\nu_f(B)+v_{k+1}$ then
$\nu_f(\omega)=e_k\nu_k(\omega)$ and consequently,
$\nu_f(\omega)\in\Lambda_f\setminus\Gamma_f$ if, and only if,
$\nu_k(\omega)\in\Lambda_k\setminus\Gamma_k$. In addition, as
$\nu_k(\omega)\leq \mu_k-1$ we have that $\nu_f(\omega)\leq
e_k(\mu_k-1)<v_{k+1}$.

If $\nu_f(H)+v_0>\nu_f(B)+v_{k+1}$, then
$\nu_k(B)-\frac{v_0}{e_k}=\nu_k(B')-\frac{v_0}{e_k}=
\nu_k(M)\in\Gamma_k$ and
$\alpha=\mu_k-1+\frac{v_0}{e_k}-\nu_k(B)\not\in\Gamma_k$.

If $\alpha\in\Lambda_k\setminus\Gamma_k$, then similarly to
Corollary \ref{valores1} we get
$v_{k+1}-e_k\alpha=\beta_{k+1}+e_k\nu_k(B)\in\Lambda_f\setminus\Gamma_f$.

On the other hand if $\alpha\not\in\Lambda_k\setminus\Gamma_k$, by
Theorem \ref{conjuntos} we have $\nu_k(M)\in\nu_k(J(f_k):f_k)$. So,
by Proposition \ref{torcao}, there exists
$\omega'\in\mathcal{E}(f_k)\cap f_k\cdot\Omega(log\ C_k)$ such that
$\nu_f(\omega')=v_{k+1}-e_k\alpha=\nu_f(\omega)$.
\end{remark}

\begin{remark}\label{remark2} As $e_k(\mu_k-1)+v_0=v_{k+1}-\beta_{k+1}$ and, by (\ref{valor}),
 we have $\nu_k(H)=\nu_k(\omega)+\mu_k-1$ it follows that
 $$\nu_f(H)+v_0<\nu_f(B)+v_{k+1}\ \ \Leftrightarrow\ \ e_k(\nu_k(\omega)-\nu_k(B))<\beta_{k+1}.$$
\end{remark}

Our aim is to define an injective function
$\rho_k:\Lambda_k\setminus\Gamma_k\rightarrow\Lambda_f\setminus\Gamma_f$.

The previous lemma give us a way to relate such sets, but it depends
on the expression of $1$-form. For instance, if
$\delta_k=\nu_k(Adx-Bdy)\in\Lambda_k\setminus\Gamma_k$ is such
$e_k(\delta_k-\nu_k(B))<\beta_{k+1}$ then we can have
$\omega_i=A_idx-B_idy\in\mathcal{E}(f_k)\setminus
f_k\cdot\Omega(log\ C_k)$ with $i=1,2$,
$\nu_k(\omega_1)=\delta_k=\nu_k(\omega_2)$ and
$\nu_k(B_1)\neq\nu_k(B_2)$.

In order to deal with this situation we consider the following
function considered by Delorme in \cite{delorme}:
$$\begin{array}{cccl}
\Theta_k : & \Lambda_k & \rightarrow & \Gamma_k\cup\{\infty\} \\
 & \delta_k & \mapsto & \max\{\nu_k(B);\ \delta_k=\nu_k(Adx-Bdy)\}.
\end{array}$$

\begin{remark}\label{not-injective}
The function $\Theta_k$ is not injective because
$\Theta_k(\delta_k)=\infty$ for every
$\delta_k=\nu_k(Adx)\in\Gamma_k+\frac{v_0}{e_k}\subset\Gamma_k$.

In addition if $\delta_k=\sum_{i=1}^{k}s_i\frac{v_i}{e_k}$, then
$\delta_k=\nu_k\left ( d\left ( \prod_{i=1}^{k}f_{i-1}^{s_i}\right
)\right )$ where $\{f_0,\ldots ,f_k\}$ is a complete system of
semiroots of $f_k$ and $0\leq s_i<n_{i}$ for $i=1,\ldots ,k$. As
$\nu_k\left ( \prod_{i=1, i\neq
j}^{k}f_{i-1}^{s_i}f_{j-1}^{s_j-1}(f_{j-1})_y\right
)=\delta_k-\frac{\beta_j}{e_k}$ for $s_j\neq 0$, we have that
$\Theta_k(\delta_k)\geq\delta_k-\frac{\beta_r}{e_k}$ where
$r=\max\{j;\ s_j\neq 0\}$.
\end{remark}

However, $\Theta_k:\Lambda_k\setminus\Gamma_k\rightarrow\Gamma_k$ is
injective. In fact, given $\delta_1,
\delta_2\in\Lambda_k\setminus\Gamma_k$ with $\delta_1<\delta_2$ if
$\Theta_k(\delta_1)=\Theta_k(\delta_2)$, then there exist
$\omega_i=A_idx-B_idy$ with $\nu_k(\omega_i)=\delta_i$ and
$\nu_k(B_i)=\Theta(\delta_i)$ for $i=1,2$. In this way, there exists
$c\in\mathbb{C}$ such that
$\nu_k(B_1-cB_2)>\nu_k(B_2)=\nu_k(B_1)=\Theta_k(\delta_1)$ and
considering $\omega=\omega_1-c\omega_2=(A_1-cA_2)dx-(B_1-cB_2)dy$ we
get $\nu_k(\omega)=\nu_k(\omega_1)=\delta_1$ and
$\Theta_k(\delta_1)\geq \nu_k(B_1-cB_2)>\Theta_k(\delta_1)$ that is
an absurd.

We define
\begin{equation}\label{rho}
\begin{array}{cccl}
\rho_k : & \Lambda_k\setminus\Gamma_k & \rightarrow &
\Lambda_f\setminus\Gamma_f \\
 & \delta_k & \mapsto & \left \{
 \begin{array}{lcl}
e_k\delta_k & \mbox{if} &
e_k(\delta_k-\Theta_k(\delta_k))<\beta_{k+1} \\
\beta_{k+1}+e_k\Theta_k(\delta_k) & \mbox{if} &
e_k(\delta_k-\Theta_k(\delta_k))>\beta_{k+1}.
 \end{array}\right .
\end{array}
\end{equation}
As $\Theta_k:\Lambda_k\setminus\Gamma_k\rightarrow\Gamma_k$ is well
defined and injective, by Remark \ref{remark}, $\rho_k$ is well
defined and injective also.

In what follows we denote
$$L_k^1=\{\delta_k\in\Lambda_k\setminus\Gamma_k;\ e_k(\delta_k-\Theta_k(\delta_k))<\beta_{k+1}\}
\ \ \mbox{and}\ \ L_k^2=\{\delta_k\in\Lambda_k\setminus\Gamma_k;\
e_k(\delta_k-\Theta_k(\delta_k))>\beta_{k+1}\}.
$$

The following proposition summarizes the discussion of this
subsection.

\begin{proposition}\label{valores2}
With the notations introduced in this subsection we have that
$\rho_k(L_k^1)=e_kL_k^1\subset\Lambda_f\setminus\Gamma_f$ and
$$\rho_k(L_k^2)=\left \{\beta_{k+1}+e_k\Theta_k(\delta_k)=v_{k+1}-
e_k\left ( \mu_k-1+\frac{v_0}{e_k}-\Theta_k(\delta_k)\right );\
\delta_k\in L_k^2\right \}\subset\Lambda_f\setminus\Gamma_f.$$ In
addition,
$\sum_{i=k+1}^{g}s_iv_i+\rho_k(L_k^2)\subset\Lambda_f\setminus\Gamma_f$
for $0\leq s_i<n_i$ for $i=k+2,\ldots ,g$ and $0\leq s_{k+1}\leq
n_{k+1}-2$.
\end{proposition}
\Dem The result follows directly form the above analysis and,
similarly the computations done in Corollary \ref{valores1}, we get
$\sum_{i=k+1}^{g}s_iv_i+\rho_k(L_k^2)\subset\Lambda_f\setminus\Gamma_f$
for $0\leq s_i<n_i$ for $i=k+2,\ldots ,g$ and $0\leq s_{k+1}\leq
n_{k+1}-2$. \cqd

\subsection{Main results}

The results presented in this section allow us to determine elements
in $\Lambda_f\setminus\Gamma_f$ by $\Lambda_k\setminus\Gamma_k$ and
consequently to relate the Tjurina number $\tau_f$ of $C_f$ with the
Tjurina number $\tau_k$ of $C_{k}$.

\begin{theorem}\label{contagem} For any branch $C_f$ such that $f_k$ is a
$k$-semiroot of $f$ with $0\leq k< g$ we have
$$e_kL_k^1\ \dot{\cup}\ \left \{ \sum_{i=k+1}^{g}s_iv_i+\rho_k(L_k^2)\right \}
\ \dot{\cup}\ \left \{
\sum_{i=k+1}^{g}s_iv_i+v_{k+1}-e_k\delta\right
\}\subseteq\Lambda_f\setminus\Gamma_f$$ with
$$0\leq s_i<n_i; \ i=k+2,\ldots, g, \ 0\leq s_{k+1}\leq n_{k+1}-2 \ \mbox{and} \ \delta\in\{\alpha, \  \alpha\in\mathbb{N}^{\ast}\setminus\Lambda_k
\}\ \dot{\cup}\ \{-\alpha, \ \alpha\in\mathbb{ N }\setminus\Gamma_k
\}.$$
\end{theorem}
\Dem By Corollary \ref{valores1} and Proposition \ref{valores2} we
have
$$\left \{
\sum_{i=k+1}^{g}s_iv_i+v_{k+1}-e_k\delta\right
\}\subseteq\Lambda_f\setminus\Gamma_f\ \ \ \ \mbox{and}\ \ \ \
e_kL_k^1\ \dot{\cup}\ \left \{
\sum_{i=k+1}^{g}s_iv_i+\rho_k(L_k^2)\right
\}\subseteq\Lambda_f\setminus\Gamma_f.$$

It is immediate that
$\left\{\displaystyle\sum_{i=k+1}^gs_iv_i+v_{k+1}-e_{k}\delta
\right\}\cap e_{k}L_k^1=\emptyset$.

If $\left \{ \sum_{i=k+1}^{g}s_iv_i+\rho_k(L_k^2)\right \} \ \cap\
\left \{ \sum_{i=k+1}^{g}s_iv_i+v_{k+1}-e_k\delta\right
\}\neq\emptyset$ then we must have
$$\sum_{i=k+1}^gs_iv_i+v_{k+1}-e_{k}\gamma=\sum_{i=k+1}^gr_iv_i+v_{k+1}-e_{k}\delta,
$$
with $v_{k+1}-e_k\gamma\in\rho_k(L_k^2)$, $0\leq r_i,s_i<n_i; \
i=k+2,\ldots, g, \ 0\leq r_{k+1},s_{k+1}\leq n_{k+1}-2$. By Remark
\ref{padrao}, we get $\gamma=\delta$.

On the other hand, by Proposition \ref{torcao} and Proposition
\ref{valores2} there exist $\omega_1=A_1dx-B_1dy\in
\mathcal{E}(f_k)\cap f_k\cdot\Omega(log\ C_k)$ and
$\omega_2=A_2dx-B_2dy\in \mathcal{E}(f_k)\setminus f_k\cdot\Omega
(log\ C_k)$ such that
$\mu_k-1+\frac{v_0}{e_k}-\nu_k(B_1)=\delta=\gamma=\mu_k-1+\frac{v_0}{e_k}-\nu_k(B_2)$.
So, $\nu_k(B_1)=\nu_k(B_2)=\Theta_k(\nu_k(\omega_2))$.

Taking $c\in\mathbb{C}$ such that
$\nu_k(B_1-cB_2)>\nu_k(B_1)=\nu_k(B_2)$ we obtain
$\omega=\omega_1-c\omega_2\in\mathcal{E}(f_k)\setminus
f_k\cdot\Omega(log\ f_k)$ with $\nu_k(\omega)=\nu_k(\omega_2)$ and
$\nu_k(B_2)=\Theta_k(\nu_k(\omega_2))=\Theta_k(\nu_k(\omega))\geq
\nu_k(B_1-cB_2)>\nu_k(B_2)$, that is an absurd.

Hence, $\left \{ \sum_{i=k+1}^{g}s_iv_i+\rho_k(L_k^2)\right \} \
\cap\ \left \{ \sum_{i=k+1}^{g}s_iv_i+v_{k+1}-e_k\delta\right
\}=\emptyset$ and we conclude the proof.
 \cqd

For any branch
$C_f$ we have that $\tau_f=\mu_f-\sharp (\Lambda_f\setminus\Gamma_f)$ and consequently we obtain the following corollary.

\begin{corollary}\label{cor-tjurina1} For any branch $C_f$ such that $f_k$ is a
$k$-semiroot of $f$ with $0\leq k< g$ we have
$\sharp(\Lambda_f\setminus\Gamma_k)\geq\mu_k+e_{k+1}(n_{k+1}-2)\tau_k$.
In particular, $\tau_f\leq\mu_f-\mu_k-e_{k+1}(n_{k+1}-2)\tau_k.$
\end{corollary}
\Dem Notice that $\Lambda_k\setminus\Gamma_k=L_k^1\ \dot{\cup}\
L_k^2$, $\sharp (\rho_k(L_k^1))=\sharp(L_k^1)$, {\small $$\sharp\left
\{ \sum_{i=k+1}^{g}s_iv_i+\rho_k(L_k^2)\right
\}=e_{k+1}(n_{k+1}-1)\sharp(L_k^2)\ \ \mbox{and}\ \ \sharp\left \{
\sum_{i=k+1}^{g}s_iv_i+v_{k+1}-e_k\delta\right
\}=e_{k+1}(n_{k+1}-1).\sharp (T)$$} with
$$0\leq s_i<n_i; \ i=k+2,\ldots, g, \ 0\leq s_{k+1}\leq n_{k+1}-2 \ \mbox{and} \ \delta\in  T:=\{\alpha, \  \alpha\in \mathbb{N}^{\ast}\setminus\Lambda_k
\}\ \dot{\cup}\ \{-\alpha, \ \alpha\in\mathbb{ N }\setminus\Gamma_k
\}.$$

As $\sharp(\mathbb{N}\setminus\Gamma_k)=\frac{\mu}{2}$ (see
(\ref{milnor})) and
$\mathbb{N}^*\setminus\Lambda_k=(\mathbb{N}\setminus
\Gamma_k)\setminus (\Lambda_k\setminus\Gamma_k)$ we get $\sharp
(T)=\mu_k-\sharp(\Lambda_k\setminus\Gamma_k)=\tau_k$.

By the previous theorem we obtain
$$\begin{array}{cl}
\sharp(\Lambda_f\setminus\Gamma_f) & \geq \sharp(L_k^1))+
e_{k+1}(n_{k+1}-1)(\sharp(L_k^2)+\tau_k) \\
& =\sharp (L_k^1)+\sharp
(L_k^2)+\tau_k+e_{k+1}(n_{k+1}-2)(\sharp(L_k^2)+\tau_k)
\\
&
=\sharp(\Lambda_k\setminus\Gamma_k)+\tau_k+e_{k+1}(n_{k+1}-2)(\sharp(L_k^2)+\tau_k)
\\
& =\mu_k+e_{k+1}(n_{k+1}-2)(\sharp(L_k^2)+\tau_k)\geq
\mu_k+e_{k+1}(n_{k+1}-2)\tau_k.\end{array}$$ In
particular, as $\tau_f=\mu_f-\sharp(\Lambda_f\setminus\Gamma_f)$ we
get $\tau_f\leq \mu_f-\mu_k-e_{k+1}(n_{k+1}-2)\tau_k.$ \cqd

As $\frac{3\mu_f}{4}\leq \tau_f$ for any branch
$C_f$ (see \cite{alberich},
\cite{genzmer} and \cite{japa}) considering the semiroot $f_{g-1}$, the previous results give
us an upper bound for $\tau_f$ in terms of topological data.

\begin{corollary}\label{corolario} For any branch $C_f$ with semigroup $\Gamma_f=\langle v_0,\ldots ,v_g\rangle$
we have $\tau_f\leq \mu_f-\frac{(3n_g-2)}{4}\mu_{g-1}$.
\end{corollary}
\Dem As $e_g=1$, by Corollary \ref{cor-tjurina1} and using that
$\frac{3\mu_{g-1}}{4}\leq \tau_{g-1}$ we obtain the result. \cqd

\begin{remark} In \cite{genzmer}, the authors proof that
$\tau_f\geq\frac{3}{4}\mu_f+\left (
\frac{-1+\sqrt{1+4\mu_f}}{8}\right ) $ for any irreducible plane
curve $C_f$. So, we can use this inequality for $f_{g-1}$ to obtain a finer topological upper bound for $\tau_f$ in the above corollary.
\end{remark}

\begin{example} In \cite{6919} the authors consider branches $C_f$ with
semigroup $\Gamma=\langle 6,9,19\rangle$ and they present all
possibilities for $\Lambda_f\setminus\Gamma_f$ and consequently all
possible values for $\tau_f$.

Any branch admitting this semigroup $\Gamma$ can be done by a
parametrization $(t^6,t^9+at^9+ \mbox{(h.o.t.)})$ with $a\neq 0$ and
all of them share the semiroot $f_1=y^2-x^3$ with $\Gamma_1=\langle
2,3\rangle$, $\Lambda_1\setminus\Gamma_1=\emptyset$ and
$\tau_1=\mu_1=2$.

Considering $k=1$ in Theorem \ref{contagem}, we obtain that
$\{16,22,35,41\}\subseteq\Lambda_f\setminus\Gamma_f$ and by
Corollary \ref{cor-tjurina1} (or Corollary \ref{corolario})
$\tau_f\leq 38$ for any branch $C_f$. In fact, in \cite{6919} we
found the following possibilites:
\begin{center}
	\begin{tabular}{|c|c|}
	\hline
	$\Lambda\setminus\Gamma$ & $\tau$ \\
	\hline	
$\{16, 22, 26, 29, 32, 35, 41\}$ & $35$ \\	
\hline
$\{16, 22, 26, 32, 35, 41\}$ & $36$ \\	
\hline
$\{16, 22, 29, 32, 35, 41\}$ & $36$ \\	
\hline
$\{16, 22, 29, 35, 41\}$ & $37$ \\	
\hline
    \end{tabular}
\end{center}
\end{example}

Notice that in above example we get $n_g=n_2=3$. In the next section
we will show that for semigroup $\Gamma=\langle v_0,\ldots
,v_g\rangle$ with $n_g=2$ we obtain the equality in Corollary
\ref{corolario}.

\section{Branches with semigroup $\Gamma =\langle v_0,\ldots ,v_g\rangle$ and $n_g=2$}

In \cite{luengo}, Luengo and Pfister proof that in the topological
class determined by semigroup $\Gamma=\langle v_0,v_1,v_2\rangle$
with $n_2=GCD(v_0,v_1)=2$ the Tjurina number $\tau_f$ is the same
for any branch $C_f$ belonging to it and $\tau_f=\mu_f-\left (
\frac{v_0}{2}-1\right )\left (\frac{v_1}{2}-1 \right
)=\mu_f-\mu_{1}$ where, as before, $\mu_1$ denote the Milnor number
of the $1$-semiroot $f_1$. Inspired by the result of Luengo and
Pfister, Watari (see \cite{watari}) ask if the Tjurina number can be
expressed using topological data for branches with semigroup
$\Gamma=\langle v_0,v_1,\ldots ,v_g\rangle$ with $n_i=2$ for
$i=1,\ldots ,g$.

In this section we will show that the inequality in Corollary
\ref{corolario} is sharp for any branch $C_f$ in the equisingularity
class determined by $\Gamma=\langle v_0,v_1,\ldots ,v_g\rangle$ with
$n_g=2$, that is, $\tau_f=\mu_f-\mu_{g-1}$. In particular, we obtain
a generalization of the result of Luengo-Pfister and an answer for a
Watari's question.

Notice that for $g=1$ if $n_1=2$ then we get $\Gamma=\langle
2,v_1\rangle$ and $\tau_f=\mu_f=v_1-1$ for any curve $C_f$ with
semigroup $\Gamma$ (see \cite{torsion} for instance). So, in what
follows we will suppose that $g>1$.

Firstly we present a result for any branch with semigroup $\Gamma=\langle v_0,\ldots ,v_g\rangle$ without restriction on $n_i$.

\begin{proposition}\label{prop} With the same notation as introduced in previous section and $1\leq k<g$ we have that
    $$\{\lambda\in\Lambda_f\setminus\Gamma_f;\ \lambda <v_{k+1}\}=\rho_k(\Lambda_k\setminus\Gamma_k)\ \dot{\cup}\ \{v_{k+1}-e_k\delta_k;\ \delta_k\in\mathbb{N}^*\setminus\Lambda_{k}\}.$$
\end{proposition}
\Dem By Corollary \ref{valores1}, using the map $\rho_k$ defined in
(\ref{rho}) and Theorem \ref{contagem} we get
$$\rho_k(\Lambda_k\setminus\Gamma_k)\ \dot{\cup}\ \{v_{k+1}-e_k\delta_k;
\
\delta_k\in\mathbb{N}^*\setminus\Lambda_{k}\}\subseteq\{\lambda\in\Lambda_f\setminus\Gamma_f;\
\lambda <v_{k+1}\}.$$

To guarantee the other inclusion we consider
$\lambda=\nu_f(\omega)\in\Lambda_f\setminus\Gamma_f$ with
$\lambda<v_{k+1}$ with $\omega\in\Omega^1$. In particular, we can
assume that $\omega\in\mathcal{E}(f_k)\setminus\mathcal{F}(f_k)$.

If $\omega\in\mathcal{E}(f_k)\cap f_k\cdot\Omega(log\ C_k)$, by
Corollary \ref{valores1}, we have that
$\nu_f(\omega)=v_{k+1}-e_k\delta_k$ with
$\delta_k\in\mathbb{N}^*\setminus\Lambda_k$.

Suppose that $\omega=Adx+Bdy\in\mathcal{E}(f_k)\setminus
f_k\cdot\Omega (log\ C_k)$.

If $e_k(\nu_k(\omega)-\nu_k(B))<\beta_{k+1}$ then, by Remark
\ref{remark}  and Remark \ref{remark2}, we obtain
$\lambda=e_k\nu_k(\omega)\in\rho_k(L_k^1)\subseteq\rho_k(\Lambda_k\setminus\Gamma_k)$.

On the other hand if $e_k(\nu_k(\omega)-\nu_k(B))>\beta_{k+1}$, by
Lemma \ref{4.3} and Remark \ref{remark2} we have
$\lambda=v_{k+1}-e_k\alpha$ with
$0<\alpha=\mu_k-1+\frac{v_0}{e_k}-\nu_k(B)\not\in\Gamma_k$.

If $\alpha\not\in\Lambda_k$ then
$\lambda=v_{k+1}-e_k\alpha\in\{v_{k+1}-e_k\delta_k;\
\delta_k\in\mathbb{N}^*\setminus\Lambda_k\}$.

The remain case to consider is
$\omega=Adx+Bdy\in\mathcal{E}(f_k)\setminus f_k\cdot\Omega (log\
C_k)$ with $e_k(\nu_k(\omega)-\nu_k(B))>\beta_{k+1}$,
$\nu_k(\omega)=\lambda=v_{k+1}-e_k\alpha$ and
$0<\alpha=\mu_k-1+\frac{v_0}{e_k}-\nu_k(B)\in\Lambda_k\setminus\Gamma_k$.
In particular, $\nu_k(B)-\frac{v_0}{e_k}\in\Gamma_k$, that is,
$\nu_k(B)\in\Gamma_k+\frac{v_0}{e_k}$.

We will show that
$\lambda\in\rho_k(L_k^2)\subset\rho_k(\Lambda_k\setminus\Gamma_k)$.

As $\nu_k(B)\in\Gamma_k+\frac{v_0}{e_k}$, if
$\nu_k(B)=\Theta_k(\nu_k(\omega))$ then
$\nu_k(\omega)\in\Lambda_k\setminus\Gamma_k$. In fact, if
$\nu_k(\omega)\in\Gamma_k$ then, by Remark \ref{not-injective}, we
get $\nu_k(B)=\Theta_k(\nu_k(\omega))=\infty$ or
$\nu_k(B)=\Theta_k(\nu_k(\omega))\geq
\nu_k(\omega)-\frac{\beta_k}{e_k}$, that is,
$e_k(\nu_k(\omega)-\nu_k(B))\leq \beta_k<\beta_{k+1}$ in both
situations we have a contradiction.

In this way, if $\nu_k(B)=\Theta_k(\nu_k(\omega))$ then
$\nu_k(\omega)\in L_k^2\subset\Lambda_k\setminus\Gamma_k$ and
consequently
$\gamma=v_{k+1}-e_k\alpha\in\rho_k(L_k^2)\subset\rho_k(\Lambda_k\setminus\Gamma_k)$.

To conclude the proof we consider $\nu_k(B)<\Theta_k(\nu_k(\omega))$
and $\omega_0=A_0dx-B_0dy\in\mathcal{E}(f_k)\setminus
f_k\cdot\Omega(log\ C_f)$ with
$\nu_k(\omega_0)=\nu_k(\omega)=\delta_0\in\Lambda_k\setminus\Gamma_k$,
$\nu_k(B_0)=\Theta_k(\delta_0)$ and
$\delta_1=\nu_k(\omega-\omega_0)>\nu_k(\omega)=\nu_k(\omega_0)$.
Notice that $\nu_k(B)=\nu_k(B-B_0)\leq\Theta_k(\delta_1)$.

If $\nu_k(B)<\Theta_k(\delta_1)$ we repeat the procedure again. We
claim that after a finite number of steps we obtain
$\tilde{\omega}=\tilde{A}dx-\tilde{B}dy\in\mathcal{E}(f_k)\setminus
f_k\cdot\Omega(log\ C_k)$ with
$\nu_k(\tilde{\omega})=\tilde{\delta}$ and
$\Theta_k(\tilde{\delta})=\nu_k(\tilde{B})=\nu_k(B)$.

In fact, if it is not the case we obtain
$$\omega-\sum_{i=0}^{\infty}\omega_i=\left ( A-\sum_{i=0}^{\infty}A_i\right )dx-\left ( B-\sum_{i=0}^{j}B_i\right )dy$$
with $\nu_k(\omega-\sum_{i=0}^{\infty}\omega_i)=\infty$ and
$\nu_k(B)=\nu_k(B-\sum_{i=0}^{\infty}B_i)$, that is,
$\omega-\sum_{i=0}^{\infty}\omega_i\in\mathcal{E}(f_k)\cap
f_k\cdot\Omega(log\ C_k)$. By Proposition \ref{torcao} we must have
$\mu_k-1+\frac{v_0}{e_k}-\nu_k(B)=\alpha\not\in\Lambda_k$ that it is
a contradiction.

In this way there exists
$\tilde{\omega}=\tilde{A}dx-\tilde{B}dy\in\mathcal{E}(f_k)\setminus
f_k\cdot\Omega(log\ C_k)$ with
$\nu_k(\tilde{\omega})=\tilde{\delta}\geq\nu_k(\omega)$ and
$\Theta_k(\tilde{\delta})=\nu_k(\tilde{B})=\nu_k(B)$. In particular,
$e_k(\tilde(\delta)-\Theta_k(\tilde{\delta}))\geq
e_k(\nu_k(\omega)-\nu_k(B))>\beta_{k+1}$ and, by Remark
\ref{not-injective},
$\nu_k(\tilde{\omega})\in\Lambda_k\setminus\Gamma_k$. So,
$\lambda=v_{k+1}-e_k\alpha=\rho_k(L_k^2)\subset\rho_k(\Lambda_k\setminus\Gamma_k)$
and we conclude the proof. \cqd

As a consequence any branch with $\Gamma=\langle v_0,\ldots
,v_g\rangle$ and $n_g=2$ has the same Tjurina number and it can be
expressed by means $\Gamma$.

\begin{theorem} For any branch $C_f$ in the topological class determined
by $\Gamma_f=\langle v_0,\ldots ,v_g\rangle$ with $n_g=2$ we have that
{\small $$\Lambda_f\setminus\Gamma_f=\rho_{g-1}(\Lambda_{g-1}\setminus\Gamma_{g-1})\ \dot{\cup}\
\{v_g-2\delta,\ \delta\in\mathbb{N}^*\setminus\Lambda_{g-1}\}\
\dot{\cup}\ \{v_g+2\delta;\
\delta\in\mathbb{N}\setminus\Gamma_{g-1}\}\ \ \mbox{and}\ \ \sharp
(\Lambda_f\setminus\Gamma_f)=\mu_{g-1}.$$} In particular, we have that
$\tau_f=\mu_f-\mu_{g-1}.$
\end{theorem}
\Dem Given $\Gamma_f=\langle v_0,\ldots ,v_g\rangle$ with
$n_g=e_{g-1}=2$ any $v_g<\lambda\not\in\Gamma_f$ can be uniquely
expressed as $\lambda=v_g+\sum_{i=1}^{g-1}s_iv_i-s_0v_0=v_g+2\left (
\sum_{i=1}^{g-1}s_i\frac{v_i}{2}-s_0\frac{v_0}{2}\right )$  with
$0\leq s_i<n_i$ for $i>0$ and $s_0\geq 0$ (see Remark \ref{padrao}).
In particular,
$\sum_{i=1}^{g-1}s_i\frac{v_i}{2}-s_0\frac{v_0}{2}\in\mathbb{N}\setminus\Gamma_{g-1}$
and any $v_g<\lambda\not\in\Gamma_f$ is such that
$\lambda\in\{v_g+2\delta;\
\delta\in\mathbb{N}\setminus\Gamma_{g-1}\}$, that is,
$\{\lambda\in\Lambda_f\setminus\Gamma_f;\
\lambda>v_g\}\subset\{v_g+2\gamma;\
\gamma\in\mathbb{N}\setminus\Gamma_{g-1}\}$. Consequently, by
(\ref{log-value}) we have $\{\lambda\in\Lambda_f\setminus\Gamma_f;\
\lambda>v_g\}=\{v_g+2\gamma;\
\gamma\in\mathbb{N}\setminus\Gamma_{g-1}\}$.

Now, by the previous proposition we get
$$\Lambda_f\setminus\Gamma_f=\rho_{g-1}(\Lambda_{g-1}\setminus\Gamma_{g-1})\ \dot{\cup}\
\{v_g-2\delta,\ \delta\in\mathbb{N}^*\setminus\Lambda_{g-1}\} \
\dot{\cup}\ \{v_g+2\gamma;\
\gamma\in\mathbb{N}\setminus\Gamma_{g-1}\}.$$

As
$\sharp\rho_{g-1}(\Lambda_{g-1}\setminus\Gamma_{g-1})=\sharp(\Lambda_{g-1}\setminus\Gamma_{g-1})$,
and $\sharp (\mathbb{N}\setminus\Gamma_{g-1})=\frac{\mu_{g-1}}{2}$,
we obtain
$$\sharp(\Lambda_f\setminus\Gamma_f)=\sharp(\Lambda_{g-1}\setminus\Gamma_{g-1})+
\sharp(\mathbb{N}^*\setminus\Lambda_{g-1})+
\sharp(\mathbb{N}\setminus\Gamma_{g-1})=2\sharp(\mathbb{N}\setminus\Gamma_{g-1})=\mu_{g-1}.
$$
In particular,
$\tau_f=\mu_f-\sharp(\Lambda_f\setminus\Gamma_f)=\mu_f-\mu_{g-1}$.
 \cqd

\vspace{0.5cm}

\begin{tabular}{lcl}
Abreu, M. O. R. & & Hernandes, M. E. \\
{\it osnar$@$outlook.com} & & {\it mehernandes$@$uem.br}\\
\end{tabular}
\vspace{0.5cm}

\hspace{1cm} Universidade Estadual de Maring\'a\vspace{0.25cm}

\hspace{1.5cm} Maring\'a - Paran\'a - Brazil

\end{document}